\documentclass[a4paper,submit]{aiaa}

\usepackage{graphicx}
\usepackage[english]{babel}
\usepackage{amsmath}
\usepackage{amssymb}
\usepackage{amsfonts}
\usepackage{ifthen}
\usepackage[latin1]{inputenc}
\usepackage{wasysym}
\usepackage{color}
\usepackage{nomencl}
\usepackage{colortbl}
\usepackage[T1]{fontenc}
\usepackage{ae,aecompl}
\usepackage[active]{srcltx} 
\usepackage{marvosym}
\usepackage{multirow}
\usepackage{hhline}
\usepackage{setspace}
\usepackage{supertabular}
\usepackage{url}
\usepackage{fancybox}
 \usepackage{indentfirst}
\usepackage[bookmarks=true,colorlinks,urlcolor=black,linkcolor=black,citecolor=black]{hyperref}
\usepackage{booktabs}

\def\venere   {\leavevmode\raise0.2ex\hbox{\wasyfamily\char25}}
\def\marte    {\leavevmode\lower0.2ex\hbox{\wasyfamily\char26}}

\newcommand{\V}[1]{\ensuremath{\boldsymbol{#1}}}                
\newcommand{\CC}[1]{\ensuremath{\boldsymbol{\hat{#1}}}}         
\newcommand{\T}{\textnormal{\tiny{T}}} 

\newcommand{\Ham}{\mathcal{H}}

\newcommand{\diff}{\textnormal{d}}

\newcommand{\unit}[1]{\mathrm{#1}}  

\renewcommand{\max}{\textnormal{max}}

\allowdisplaybreaks[4]

\SubmitName{Wu and Jiang}

\title{\textsf{An Atlas of Optimal Low-Thrust Rephasing Solutions in Circular Orbit}}

\author{Di Wu%
\thanks{Postdoctoral Researcher, wu-d{\scriptsize \MVAt}tsinghua.edu.cn.}, Fanghua Jiang%
\thanks{Associate Professor, jiangfh{\scriptsize \MVAt}tsinghua.edu.cn, Senior Member AIAA (\textbf{corresponding author}).}, 
and Hexi Baoyin%
\thanks{Professor, baoyin{\scriptsize \MVAt}tsinghua.edu.cn, Senior Member AIAA.} \\
{\normalsize\itshape  School of Aerospace Engineering, Tsinghua University, 100084 Beijing,} \\
{\normalsize\itshape People's Republic of China} \\}

\begin{document}

\maketitle 

\begin{spacing}{2.0}

\section{Abstract}

In this paper, the time- and propellant-optimal low-thrust rephasing problems in circular orbit are studied to depict their solution spaces in an atlas. The number of key parameters that settle the rephasing problems is reduced by developing a set of linearized equations of motion based on the Sundman transformation and by formulating two reduced shooting functions using the minimum principle and symmetry properties. Only one key parameter is identified for the time-optimal problem, while two key parameters are obtained for the propellant-optimal one. Numerical investigation of the relationships between these parameters and shooting variables reveals that they can be depicted by some curve (or contour) maps and approximated by piecewise functions (or linear interpolations). For the relatively short- or long-term rephasing cases, some analytical time- and propellant-optimal solutions are proposed and consistent with the numerical solutions. Numerical results demonstrate that the proposed solutions can provide good initial guesses to solve the low-thrust rephasing problems with nonlinear dynamics. Moreover, the approximations of the performance indexes can be used in the preliminary mission design.

\section{I. Introduction}

A category of trajectory optimization problems, known as the orbital rephasing~\cite{gonzalo2017optimal} or station change~\cite{zhao2016initial}, has gained great attention for its potential applications in the active debris removal~\cite{izzo2018kessler}, constellation reconfiguration~\cite{palmer2007reachability}, and relocation of the geostationary satellite~\cite{kim2015constrained}. The satellite driven by the impulsive chemical or low-thrust electric propulsion system~\cite{luo2007optimal,zuiani2012direct,englander2017automated} is usually considered to rendezvous with a target in the same orbit but at a different angular position using minimum time of flight or propellant consumption. Compared with traditional chemical propulsion, the electric propulsion has the advantages of lower propellant consumption and longer acting time. It has been verified in some practical missions such as the Starlink constellation and Shijian-20 geostationary satellite~\cite{yeo2021miniaturization}. As the number of these missions increases, some efficient preliminary design and optimization techniques are needed to autonomously solve the optimal low-thrust rephasing problem~\cite{gonzalo2017optimal}. However, the optimal solutions are more complicated than the bi-impulsive solution to Lambert's problem due to the continuous thrust~\cite{jiang2017improving}. The time-optimal solution is relatively easier to seek, while the propellant-optimal solution is much challenging due to the well-known bang-bang control~\cite{taheri2018generic,Warm2021Di}, whose structure (e.g., the number of burning arcs) is strongly nonlinear with respect to the boundary constraints. In addition, the states and associated costates usually need to be numerically integrated because of the nonlinear dynamics governed by the optimal low-thrust controls, which aggravates the computational time of numerical optimization. Thus, the time- and propellant-optimal low-thrust rephasing solution spaces are investigated in this study in order to obtain fast solutions.

In literature, the optimal low-thrust rephasing problem is well-established as a nonlinear programming problem by a direct method~\cite{kim2015constrained} or a two-point boundary value problem by an indirect method~\cite{zeng2017searching,bassetto2021spiral}. Based on the numerical solutions, some interesting and practical properties were presented. Thorne et al.~\cite{thorne1997minimum} showed that the time-optimal rephasing problems have similar solutions as long as the ratios of angular position change to thrust magnitude are equal. This property supported the employment of the linearized equations of motion~\cite{palmer2007reachability,bevilacqua2014analytical}, where the ratio is identified as a key independent parameter. Recently, Gonzalo and Bombardelli~\cite{gonzalo2017optimal} developed two analytical solutions to the short- and long-term time-optimal rephasing problems in a circular orbit, using a curvilinear relative motion formulation. These two solutions provide an accurate estimation of the optimal time of flight, but the transition zone between the short-term and long-term cases was still challenging. In this work, the analytical solutions will be extended to include the transition zone so as to complete a full atlas of the time-optimal rephasing solutions.

By comparison, the low-thrust propellant-optimal problems have been solved by some numerical optimization methods~\cite{meng2019low,taheri2021costate} to obtain the optimal bang-bang control. Most of the literature on this subject focuses on improving the solution efficiency and optimality, while the properties of the solutions with different boundary constraints are seldom studied to the best of the authors' knowledge. Pontani~\cite{pontani2015symmetry,pontani2015minimum} presented some symmetry properties of the optimal linear relative trajectories, which can be used to reduce the number of the unknown optimization variables (e.g., the discrete state and control variables or the costate variables). These symmetry properties should hold according to the geometrical symmetry of the rephasing problem in two-body dynamics. Thus, the rephasing problem identified as the same-orbit rendezvous~\cite{gonzalo2017optimal} can be solved simpler than the general rendezvous problem between different orbits~\cite{wu2021minimum}. Besides, the initial guesses of costates concerned with the indirect methods can be estimated by the analytical solutions~\cite{zhao2016initial,casalino2014approximate,di2022analytical} or heuristic optimization~\cite{pontani2015minimum}. The analytical solutions are usually complicated but inaccurate expressions, while the heuristic optimization needs more computational time~\cite{zhao2016initial,pontani2015minimum}. In this work, the initial guesses of costates and propellant consumption for the rephasing problem can be simply estimated by analyzing the properties of the optimal rephasing results.

The main contribution of this work is providing an atlas of the time- and propellant-optimal low-thrust rephasing solutions in the circular orbit to depict the properties of the solution spaces. Inspired by the symmetry properties~\cite{pontani2015minimum}, we simplify the optimal control problems and establish them as two reduced two-dimensional shooting functions. Meanwhile, a set of linearized equations of motion is firstly proposed based on the Sundman transformation~\cite{sundman1913memoire} and provides a good approximation to the nonlinear dynamics for the low-thrust rephasing problem. Based on this dynamical formulation, the Euler-Lagrange equations can be analytically integrated, and the number of differential equations left in the shooting function is reduced to two. Then, the numbers of key parameters for the time- and propellant-optimal problems are identified as one and two, respectively. The solution spaces are explored by a traversal method and depicted by some curve and contour maps. Two piecewise functions are fitted to approximate the initial costates and performance indexes for the time-optimal problem, while the linear interpolation method is employed for the propellant-optimal problem. In addition, the analytical solutions are discussed for the relatively short- and long-term rephasing problems and are in good agreement with the numerical solutions. Finally, several examples show that the proposed solutions can provide good initial guesses and performance indexes for the low-thrust rephasing solutions with the nonlinear dynamics.

The rest of this paper is organized as follows. First, the optimal low-thrust rephasing problem with a general performance index is formulated in Sec.~\uppercase\expandafter{\romannumeral2}. A set of linearized equations of motion is derived by the Sundman transformation and the scaling technique. Next, the time-optimal rephasing solutions are presented in Sec.~\uppercase\expandafter{\romannumeral3}. A two-dimensional shooting function with analytical Jacobian is formulated, and the time-optimal solutions are expressed by some curve maps. Then, the propellant-optimal rephasing solutions are presented in Sec.~\uppercase\expandafter{\romannumeral4}. A similar two-dimensional shooting function is derived, while the solutions are depicted by some contour maps. Finally, the numerical examples are shown in Sec.~\uppercase\expandafter{\romannumeral5} to compare the proposed solutions with those obtained using the nonlinear dynamics, and Section~\uppercase\expandafter{\romannumeral6} concludes this paper.


\section{II. Optimal Low-Thrust Rephasing Problem Formulation}

\subsection{A. Equations of Motion}

In this work, the rephasing problem in a circular orbit is considered. The satellite is assumed to  rendezvous with a target in the same orbit but at a different angular position (i.e., a phase difference of $\Delta \theta$). In a central gravitational field that neglects the perturbation terms, the optimal rephasing trajectory of the satellite must remain in the initial orbit plane. Therefore, the first three equinoctial orbital elements $\V{x} = \left[p, \, f, \, g\right]^{\T}$ and the true longitude $L$ are used to describe the motion of the satellite, while the other two elements are set to $h = k = 0$. The equations of motion are given by~\cite{Rapid2021Di}
\begin{equation} \label{dyn}
\left\{
\begin{aligned}
\dot{\V{x}} &= \V{B}\left(\V{x}, \, L\right) \, \V{a}\\[0.3cm]
\dot{L} &= A\left(\V{x}, \, L\right)
\end{aligned}\right.
\end{equation}
where $\V{a}$ denotes the thrust acceleration. The maximum magnitude of the thrust acceleration $a_{\max}$ is assumed to be constant since the mass of the satellite remains almost unchanged during the rephasing~\cite{gonzalo2017optimal}. The expressions of ${A}$ and $3\times2$ matrix $\V{B}$ are formulated as
\begin{equation} \label{AandB}
A = w^2\sqrt{\frac{\mu}{p^3}}, \quad \V{B} = \sqrt{\frac{p}{\mu}}\left[\begin{aligned} & \quad\;\, 0 \quad\;\, \qquad\quad \;\;\,2\,p\;\;\; \quad \\[0.2cm] & \;\; \sin L \;\; \quad \frac{\left(w+1\right)\cos L + f}{w} \\[0.2cm] & -\cos L \quad \frac{\left(w+1\right)\sin L + g}{w} \\[0.2cm] \end{aligned}\right]
\end{equation}
where $w = 1+ f \cos L + g \sin L$, and $\mu$ is the gravitational parameter. These equations are nonlinear and make the optimal control problem difficult to solve. A set of linearized and scaled equations of motion based on the Sundman transformation~\cite{wu2021minimum} will be derived to obtain the approximate optimal solutions.

By introducing the Sundman transformation $\diff \, L = A\left(\V{x}, \, L\right) \diff \, t$ and changing the independent variable from the time $t$ to the true longitude $L$, the equations of motion can be transformed into
\begin{equation} \label{dyn2}
\left\{
\begin{aligned}
{\V{x}}^{\prime} &= \V{B}\left(\V{x}, \, L\right) \, \V{a} \,/\, A\left(\V{x}, \, L\right)\\[0.3cm]
{t}^{\prime} &= 1 \,/\, A\left(\V{x}, \, L\right)
\end{aligned}\right.
\end{equation}
where $\left( * \right)^{\prime}$ denotes the derivative to the true longitude $L$. Then, under the assumption that the rephasing trajectory is close to the initial circular orbit, the equations of motion can be linearized around the initial orbital elements $\CC{x} = \left[\hat{p}, \, 0, \, 0\right]^{\T}$. The length and time units are scaled such that the initial semi-latus rectum (or orbital radius) $\hat{p}$ and the gravitational parameter $\mu$ are both unities and the orbit with radius $\hat{p}$ has a period of $2 \, \pi$. A set of linearized and scaled equations of motion is derived as
\begin{equation} \label{dyn3}
\left\{
\begin{aligned}
{\Delta \V{x}}^{\prime} &= \V{B}\left(\CC{x}, \, L\right) \, \V{a} \,/\, \hat{A} \\[0.3cm]
{\Delta t}^{\prime} &= -\frac{1}{\hat{A}^2} \left(\left.\frac{\partial A}{\partial \V{x}}\right|_{\V{x} = \CC{x}}\right)^{\T} {\Delta \V{x}}
\end{aligned}\right.
\end{equation}
where ${\Delta \V{x}} = \V{x} - \CC{x}$ and $\Delta t = t - \hat{t}$. The nominal time $\hat{t}$ is calculated by $\hat{t} = t_0 + \left(L-L_0\right) / \hat{A}$, where $t_0$ and $L_0$ are the initial time and true longitude, respectively. According to Eq.~\eqref{AandB}, the quantity of $\hat{A}$ is $\hat{A} = \sqrt{\mu \,/\, \hat{p}^3} \equiv 1$, and the expression of $\V{B}\left(\CC{x}, \, L\right)$ is
\begin{equation} \label{Bnorminal}
\V{B}\left(\CC{x}, \, L\right) = \left[\begin{aligned} & \quad\;\, 0 \quad\;\, \quad \;\;\;2 \quad \\[0.2cm] & \;\; \sin L \;\; \quad 2\,\cos L \\[0.2cm] & -\cos L \quad 2\, \sin L \\[0.2cm] \end{aligned}\right]
\end{equation}

Substituting Eqs.~\eqref{AandB} and~\eqref{Bnorminal} into Eq.~\eqref{dyn3}, the equations of motion can be further simplified to
\begin{equation} \label{dyn4}
\left\{
\begin{aligned}
\Delta p^{\prime} &= 2 \, a_{\theta} \\[0.3cm]
\Delta f^{\prime} &= a_{r} \, \sin L + 2\, a_{\theta} \, \cos L \\[0.3cm]
\Delta g^{\prime} &= -a_{r} \, \cos L + 2\, a_{\theta} \, \sin L \\[0.3cm]
{\Delta t}^{\prime} &= 1.5 \, {\Delta p} - 2\, \Delta f \,\cos L - 2\,\Delta g \,\sin L 
\end{aligned}\right.
\end{equation}
where $a_{r} = a \, \sin \gamma$ and $a_{\theta} = a \, \cos \gamma$ are the components of the thrust acceleration $\V{a}$ in the local vertical/local horizontal (LVLH) frame \cite{WANG2016389Analysis}, $a = \left\|\V{a}\right\|$ is the thrust acceleration magnitude, and $\gamma$ is the orientation with respect to the transversal direction. In literature, the relative motion of the rephasing is usually formulated based on the Clohessy-Wiltshire (CW) equations~\cite{bevilacqua2014analytical} or in the curvilinear coordinates~\cite{bombardelli2017approximate,gonzalo2017optimal}. The CW equations are only applicable for the rephasing with small phase difference, while the linearized equations in the curvilinear coordinates are developed for more general cases. In this work, the linearization of Eq.~\eqref{dyn2} depends only on the assumption that the elements $\V{x}$ are close to the initial elements $\CC{x}$, and the linearized equations~\eqref{dyn4} are therefore applicable for the cases with large phase difference. This characteristic of linearized equations results from the use of the Sundman transformation. The comparison between the nonlinear equations~\eqref{dyn} and linearized equations~\eqref{dyn4} will be numerically investigated in Sec.~\uppercase\expandafter{\romannumeral5}.

\subsection{B. Optimal Control Problems}

A general performance index for the optimal rephasing problem is formulated as
\begin{equation} \label{per}
J = \int_{t_0}^{t_f} \phi\left(\V{x}, \, L, \, \V{a}, \, t\right) \diff \, t = \int_{L_0}^{L_f} \frac{\phi\left(\V{x}, \, L, \, \V{a}, \, t\right)}{A\left(\V{x}, \, L\right)} \diff \, L
\end{equation}
where $\phi$ is the running cost function to be designed, and $A = \hat{A} \equiv 1$. It holds $\phi = 1$ to minimize the total time of flight and takes $\phi = a$ to minimize the total velocity increment (or propellant consumption). The techniques for solving time-optimal and propellant-optimal control problems will be introduced in Sec.~\uppercase\expandafter{\romannumeral3} and~\uppercase\expandafter{\romannumeral4}, respectively. The final time $t_f$ is free for the time-optimal problem, while it is fixed for the propellant-optimal one. The initial times in the two problems are both fixed and set to zero without loss of generality. The other boundary constraints are
\begin{equation} \label{boundt}
\begin{aligned}
&\V{x}\left(t_0\right) = \CC{x}, \quad L\left(t_0\right) = L_0 \\[0.3cm]
&\V{x}\left(t_f\right) = \CC{x}, \quad L\left(t_f\right) = L_f
\end{aligned}
\end{equation}
where $L_f$ is the final true longitude of the target at the final time. Then, the states $\Delta \V{x}$ and $\Delta t$ should satisfy
\begin{equation} \label{boundL}
\begin{aligned}
&\Delta \V{x}\left(L_0\right) = \V{0}, \quad \Delta t\left(L_0\right) = 0 \\[0.3cm]
&\Delta \V{x}\left(L_f\right) = \V{0}, \quad \Delta t\left(L_f\right) = {\Delta t}_f
\end{aligned}
\end{equation}
where the time difference $\Delta t_f = t_f - t_0 -\left(L_f - L_0\right) = \Delta \theta \in \left[-\pi, \, \pi\right]$ denotes the phase difference. If $\Delta t_f = \Delta \theta < 0$, the initial true longitude of the target is larger than that of the satellite. The rephasing problem is now formulated as a standard optimal control problem with performance index~\eqref{per}, dynamics~\eqref{dyn4}, and boundary constraints~\eqref{boundL}.

Based on the boundary constraints~\eqref{boundL}, the integration of $\Delta t$ can be simplified to
\begin{equation} \label{int}
\begin{aligned}
\Delta t\left(L_f\right) &= \int_{L_0}^{L_f} \left(1.5 \, {\Delta p} - 2\, \Delta f \,\cos L - 2\,\Delta g \,\sin L\right) \diff \, L \\[0.2cm]
&= -\int_{L_0}^{L_f} \left(1.5 \,L\, {\Delta p}^{\prime} - 2\, \Delta f^{\prime} \,\sin L + 2\,\Delta g^{\prime} \,\cos L\right) \diff \, L\\[0.2cm]
&= \int_{L_0}^{L_f} \left(2 \, a_{r} - 3 \,L\, a_{\theta}\right) \diff \, L\end{aligned}
\end{equation}
Thus, the final time difference $\Delta t\left(L_f\right)$ can be computed by integrating a differential equation $\Delta t^{\prime} = 2\,a_{r} - 3\,L\,a_{\theta}$. Note that this equation is not equivalent to Eq.~\eqref{dyn4} and can only be used for the integration. For the circular orbit, the initial true longitude can be set to an arbitrary value by intentionally adjusting the reference frame. In the following derivation, the initial and final true longitudes are set to $L_0 = -\Delta L \,/\, 2$ and $L_f = \Delta L \,/\, 2$, respectively, where $\Delta L = L_f - L_0$ denotes the true longitude difference.

\section{III. Time-optimal Rephasing Solutions}

In this section, the solutions to the time-optimal rephasing problem are presented. Based on the minimum principle~\cite{pontryagin1987mathematical}, the optimal control problem will be transformed into a two-point boundary value problem with a reduced two-dimensional shooting function. The two shooting variables are related to only one key parameter, and the relationships between them can be approximated by two piecewise functions. Besides, the analytical solutions for the cases with relatively small and large true longitude differences $\Delta L$ will be derived. 

\subsection{A. Reduced Shooting Function}

The Hamiltonian of the time-optimal control problem is established as
\begin{equation} \label{Ham1}
\begin{aligned}
\Ham &= a_{r} \left(\lambda_{\Delta f} \sin L - \lambda_{\Delta g} \cos L\right) + 2 \, a_{\theta} \left(\lambda_{\Delta p} + \lambda_{\Delta f} \cos L + \lambda_{\Delta g} \sin L\right) \\[0.2cm]
&+1.5 \, \lambda_{\Delta t} \, \Delta p - 2\,\lambda_{\Delta t} \, \Delta f\, \cos L - 2\, \lambda_{\Delta t} \, \Delta g \, \sin L+ 1\end{aligned}
\end{equation}
where $\V{\lambda} = \left[\lambda_{\Delta p}, \, \lambda_{\Delta f}, \, \lambda_{\Delta g}, \, \lambda_{\Delta t} \right]^{\T}$ is a set of costate variables associated with the states $\Delta \V{x}$ and $\Delta t$. Based on the minimum principle, all the first-order optimality necessary conditions can be derived. The Euler-Lagrange equations are given by
\begin{equation} \label{Eul1}
\left\{
\begin{aligned}
&\lambda_{\Delta p}^{\prime} = -\frac{\partial \Ham}{\partial \Delta p} = -1.5 \lambda_{\Delta t} \\[0.3cm]
&\lambda_{\Delta f}^{\prime} = -\frac{\partial \Ham}{\partial \Delta f} = 2\,\lambda_{\Delta t} \, \cos L \\[0.3cm]
&\lambda_{\Delta g}^{\prime} = -\frac{\partial \Ham}{\partial \Delta g} = 2\,\lambda_{\Delta t} \, \sin L \\[0.3cm]
&\lambda_{\Delta t}^{\prime} = -\frac{\partial \Ham}{\partial \Delta t} = 0 \end{aligned}\right.
\end{equation}
Thus, the costate $\lambda_{\Delta t}\left(L\right) \equiv \lambda_0$ remains unchanged but is unknown. According to Ref.~\cite{pontani2015symmetry}, the optimal control and costate variables for this problem should be symmetric. At the mid-point of the optimal trajectory, the transversal thrust acceleration $a_{\theta}$ and the radial thrust acceleration derivative $a_{r}^{\prime}$ both be zero. The costates at the mid-point take the values of
\begin{equation} \label{Sym1}
\begin{aligned}
&\left.\left(\lambda_{\Delta f} \sin L - \lambda_{\Delta g} \cos L\right)^{\prime}\right|_{L = 0} = \lambda_{\Delta f}\left(0\right) = 0 \\[0.3cm]
&\left.\left(\lambda_{\Delta p} + \lambda_{\Delta f} \cos L + \lambda_{\Delta g} \sin L\right)\right|_{L = 0} = \lambda_{\Delta p}\left(0\right) = 0\end{aligned}
\end{equation}
Substituting Eq.~\eqref{Sym1} into Eq.~\eqref{Eul1} and replacing $\lambda_{\Delta t}$ by the constant $\lambda_0$, the Euler-Lagrange equations can be integrated as
\begin{equation} \label{Eul1_Solu}
\left\{
\begin{aligned}
&\lambda_{\Delta p}\left(L\right) = -1.5 \, \lambda_0 \, L \\[0.3cm]
&\lambda_{\Delta f}\left(L\right) = 2\,\lambda_0\,\sin L \\[0.3cm]
&\lambda_{\Delta g} \left(L\right) = \lambda_0\,\left(\lambda_1 \, - 2\,\cos L\right) \end{aligned}\right.
\end{equation}
where $\lambda_1$ is another unknown constant. To minimize the Hamiltonian, the optimal control should follow
\begin{equation} \label{optC}
\begin{aligned}
\quad a_{r}^{\star}\left(L\right) &= \frac{a_{\max}  \, \left(\lambda_1 \cos L - 2\right)\, \textnormal{sign} \left(\lambda_0\right)}{\sqrt{\left(3\, L - 2\,\lambda_1 \sin L\right)^2 + \left(\lambda_1 \cos L - 2\right)^2}} \\[0.3cm]
a_{\theta}^{\star}\left(L\right) &= \frac{a_{\max}  \, \left(3\,L - 2\,\lambda_1 \sin L\right)\, \textnormal{sign} \left(\lambda_0\right)}{\sqrt{\left(3\, L - 2\,\lambda_1 \sin L\right)^2 + \left(\lambda_1 \cos L - 2\right)^2}}
\end{aligned}
\end{equation}
where $\left(*\right)^{\star}$ denotes the corresponding optimal value. Based on Eq.~\eqref{optC}, the symmetry property in the optimal control can be easily obtained as $a_{r}\left(L\right) = a_{r}\left(-L\right)$ and $a_{\theta}\left(L\right) = -a_{\theta}\left(-L\right)$. The constant $\lambda_0$ affects only the signs of the thrust acceleration components, and $\lambda_1$ determines their magnitude. Since the final states are fixed, the final costates are free. Besides, the transversality condition is derived as
\begin{equation} \label{trans}
\Ham\left(L_f\right) = 1 - a_{\max} \left|\lambda_0\right| \sqrt{\left(3\, L_f - 2\,\lambda_1 \sin L_f\right)^2 + \left(\lambda_1 \cos L_f - 2\right)^2} = 0
\end{equation}
The magnitude of constant $\lambda_0$ should be
\begin{equation} \label{lambda0}
\left|\lambda_0\right| = \frac{1}{a_{\max}\sqrt{\left(3\, L_f - 2\,\lambda_1 \sin L_f\right)^2 + \left(\lambda_1 \cos L_f - 2\right)^2}}
\end{equation}
which does not actually affect the solution to the time-optimal control problem because of Eq.~\eqref{optC}. Therefore, all the first-order optimality conditions are included in the optimal control expressed by Eq.~\eqref{optC}.

The two-point boundary value problem is then formulated. After setting the sign of $\lambda_0$, the value of $\lambda_1$, and the true longitude difference $\Delta L$, the optimal control can be calculated by Eq.~\eqref{optC}, and the states can be obtained through integrating the four-dimensional ordinary differential equations~\eqref{dyn4}. The final states should satisfy the boundary constraints~\eqref{boundL}. Among them, $\Delta p \left(L_f\right) = 0$ and $\Delta f \left(L_f\right) = 0$ are automatically guaranteed by the symmetry property. Bearing in mind Eq.~\eqref{int}, the remaining constraints yield:
\begin{equation} \label{bound1}
\begin{aligned}
\Delta g \left(L_f\right) &= \textnormal{sign} \left(\lambda_0\right) a_{\max} \int_{L_0}^{L_f} \frac{6\,L\,\sin L + 2\, \cos L - \lambda_1 - 3\,\lambda_1 \, \sin^2 L}{\sqrt{\left(3\, L - 2\,\lambda_1 \sin L\right)^2 + \left(\lambda_1 \cos L - 2\right)^2}} \diff \, L = 0 \\[0.3cm]
\Delta t \left(L_f\right) &= \textnormal{sign} \left(\lambda_0\right) a_{\max} \int_{L_0}^{L_f} \frac{6\,\lambda_1\,L\,\sin L + 2\,\lambda_1\,\cos L - 9\,L^2 - 4}{\sqrt{\left(3\, L - 2\,\lambda_1 \sin L\right)^2 + \left(\lambda_1 \cos L - 2\right)^2}} \diff \, L = \Delta t_f
\end{aligned}
\end{equation}
The dimension of the ordinary differential equations is therefore reduced from four to two.

Finally, the true longitude difference $\Delta L$ and costate $\lambda_1$ are identified as the shooting variables. According to Eq.~\eqref{bound1}, a two-dimensional shooting function is obtained:
\begin{equation} \label{shooting1}
\V{\varPhi} \left(\V{z}\right) = \left[\begin{aligned}
& \;\;\, \int_{L_0}^{L_f} \frac{6\,L\,\sin L + 2\, \cos L - \lambda_1 - 3\,\lambda_1 \, \sin^2 L}{\sqrt{\left(3\, L - 2\,\lambda_1 \sin L\right)^2 + \left(\lambda_1 \cos L - 2\right)^2}} \diff \, L \\[0.2cm]
& \int_{L_0}^{L_f} \frac{ 9\,L^2 + 4 - 6\,\lambda_1\,L\,\sin L - 2\,\lambda_1\,\cos L}{\sqrt{\left(3\, L - 2\,\lambda_1 \sin L\right)^2 + \left(\lambda_1 \cos L - 2\right)^2}} \diff \, L - \chi \end{aligned}\right] = \V{0}
\end{equation}
where the key parameter $\chi = - \, \textnormal{sign} \left(\lambda_0\right) \Delta t_f \,/\, a_{\max}$ characterizes the time-optimal solution and determines the shooting variables $\V{z} = \left[\Delta L, \, \lambda_1\right]^{\T}$. Two branches of solutions, $\left\{\Delta L, \, \lambda_1, \lambda_0, \Delta t\right\}$ and $\left\{\Delta L, \, \lambda_1, -\lambda_0, -\Delta t\right\}$, can be obtained for the same parameter $\chi$. Numerical simulation shows that this parameter is always positive, and $\lambda_0$ should be negative when $\Delta t>0$. In general, this shooting function is nonlinear and difficult to solve. Since $\chi$ is the only key parameter in shooting function~\eqref{shooting1}, the numerical solutions and analytical approximations with different values of $\chi$ can be investigated by a single-variable traversal method in the next subsection.

\subsection{B. Numerical Solutions and Analytical Approximations}

To find the relationships between $\chi$, $\Delta L$, and $\lambda_1$, a traversal method is used in this work. Based on the fact that the magnitude of $\chi$ must be a monotonic increasing function with respect to the true longitude difference $\Delta L$, $\Delta L$ is chosen to be the independent variable instead of the parameter $\chi$. For each specific true longitude difference, the costate $\lambda_1$ can be obtained easily by solving the first shooting function, and the parameter $\chi$ is obtained by the second shooting function.

\begin{figure}[ht!]
	\centering
	\includegraphics[scale = 0.6]{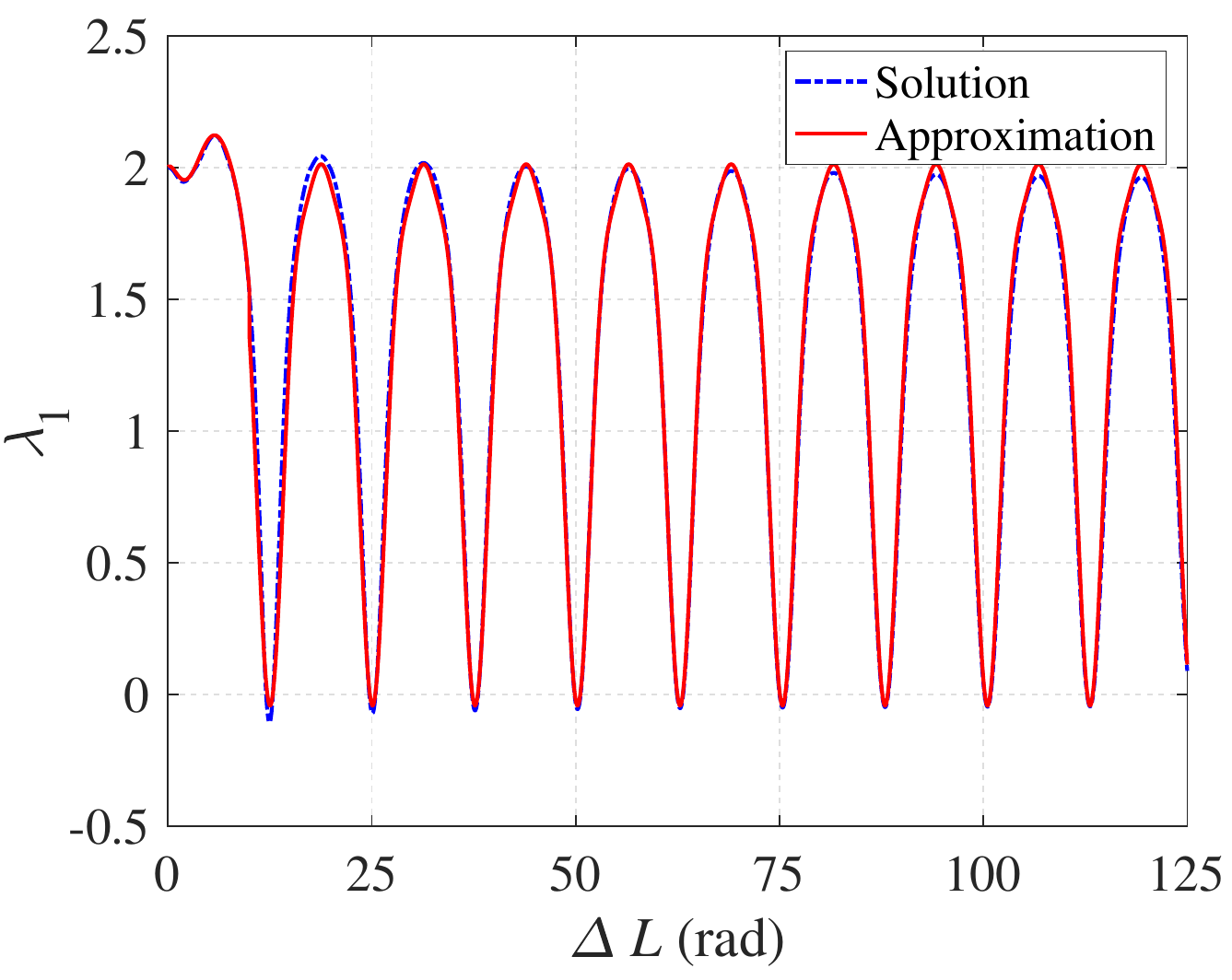}
	\caption{The curves of numerical solution and analytical approximation of $\lambda_{1}$.}
	\label{fig1}
\end{figure}

According to the symmetry property, the first shooting function can be transformed into
\begin{equation} \label{numerical1}
f_1\left(\Delta L, \, \lambda_1\right) = \int_{0}^{\Delta L/2} \frac{6\,L\,\sin L + 2\, \cos L - \lambda_1 - 3\,\lambda_1 \, \sin^2 L}{\sqrt{\left(3\, L - 2\,\lambda_1 \sin L\right)^2 + \left(\lambda_1 \cos L - 2\right)^2}} \diff \, L = 0
\end{equation}
where the true longitude difference $\Delta L$ is fixed. The partial derivatives of $f_1\left(\Delta L, \lambda_1\right)$ with respect to $\Delta L$ and $\lambda_1$ are derived as
\begin{equation} \label{partial1}
\begin{aligned}
&\frac{\partial \, f_1}{\partial \, \Delta L} = \frac{3\,\Delta L\, \sin \left(\Delta L/2\right) + 2\,\cos\left(\Delta L/2\right) - \lambda_1 - 3\,\lambda_1\,\sin^2\left(\Delta L/2\right)}{\sqrt{\left[3\,\Delta L - 4\,\lambda_1 \sin \left(\Delta L/2\right)\right]^2 + 4 \left[\lambda_1 \cos \left(\Delta L/2\right) - 2\right]^2}}\\[0.3cm]
&\frac{\partial \, f_1}{\partial \, \lambda_1} = \int_{0}^{\Delta L/2} \frac{ - \left(3 \, L \, \cos L - 4\,\sin L\right)^2}{\sqrt{\left(3\, L - 2\,\lambda_1 \sin L\right)^2 + \left(\lambda_1 \cos L - 2\right)^2}^ 3} \diff \, L < 0
\end{aligned}
\end{equation}
Thus, $f_1$ is a monotonic decreasing function with respect to $\lambda_1$. The nonlinear equation~\eqref{numerical1} of $\lambda_1$ must have a unique solution. It can be solved by a Newton's iteration method with arbitrary initial value. By traversing $\Delta L$ between 0.0125 and 125.0 in a step of 0.0125, the first shooting function~\eqref{numerical1} for each case is solved, and the numerical relationship between $\lambda_1$ and $\Delta L$ is then depicted in Fig.~\ref{fig1}. The costate $\lambda_1$ is about 2.0 when $\Delta L$ is small, and it is roughly a periodic function between 0.0 and 2.0 when $\Delta L$ is relatively large. To provide a more accurate approximation, a piecewise function based on the finite Fourier series is fitted using the MATLAB curve fitting toolbox program {\it cftool}. The analytical costate values calculated by the piecewise function are shown in Fig.~\ref{fig1} for comparison. The piecewise function is given as
\begin{equation} \label{lambda_g1}
\lambda_1 \approx \left\{\begin{aligned} &c_{10} + \sum\limits_{i=1}^{3} \left[c_{1i} \, \cos \left( i \, n_1\, \Delta L\right) + d_{1i} \, \sin \left(i\,n_1\,\Delta L\right)\right] \quad  \textnormal{if} \quad \Delta L \leq 10 \\[0.2cm]
&c_{20} + \sum\limits_{i=1}^{3} \left[c_{2i} \, \cos \left( i \, n_2\, \Delta L\right) + d_{2i} \, \sin \left(i\,n_2\,\Delta L\right)\right] \quad  \textnormal{if} \quad \Delta L > 10\end{aligned}\right.
\end{equation}
where the coefficients are obtained as $c_{10} = -19.34, \, c_{11} = 22.5, \, c_{12} = 1.261, c_{13} = -2.419, \, d_{11} = 23.9, \, d_{12} =  - 14.18, \, d_{13} = 1.54, \, n_1 = 0.1699$, $c_{20} = 1.302, \, c_{21} = - 0.9269, \, c_{22} = - 0.3164, c_{23} = - 0.09964, \, d_{21} = 0.02194, \, d_{22} = 0.01196, \, d_{23} = 0.005974$, and $n_2 = 0.4999$, respectively.

%
%

\begin{figure}[ht!]
	\centering
	\includegraphics[scale = 0.6]{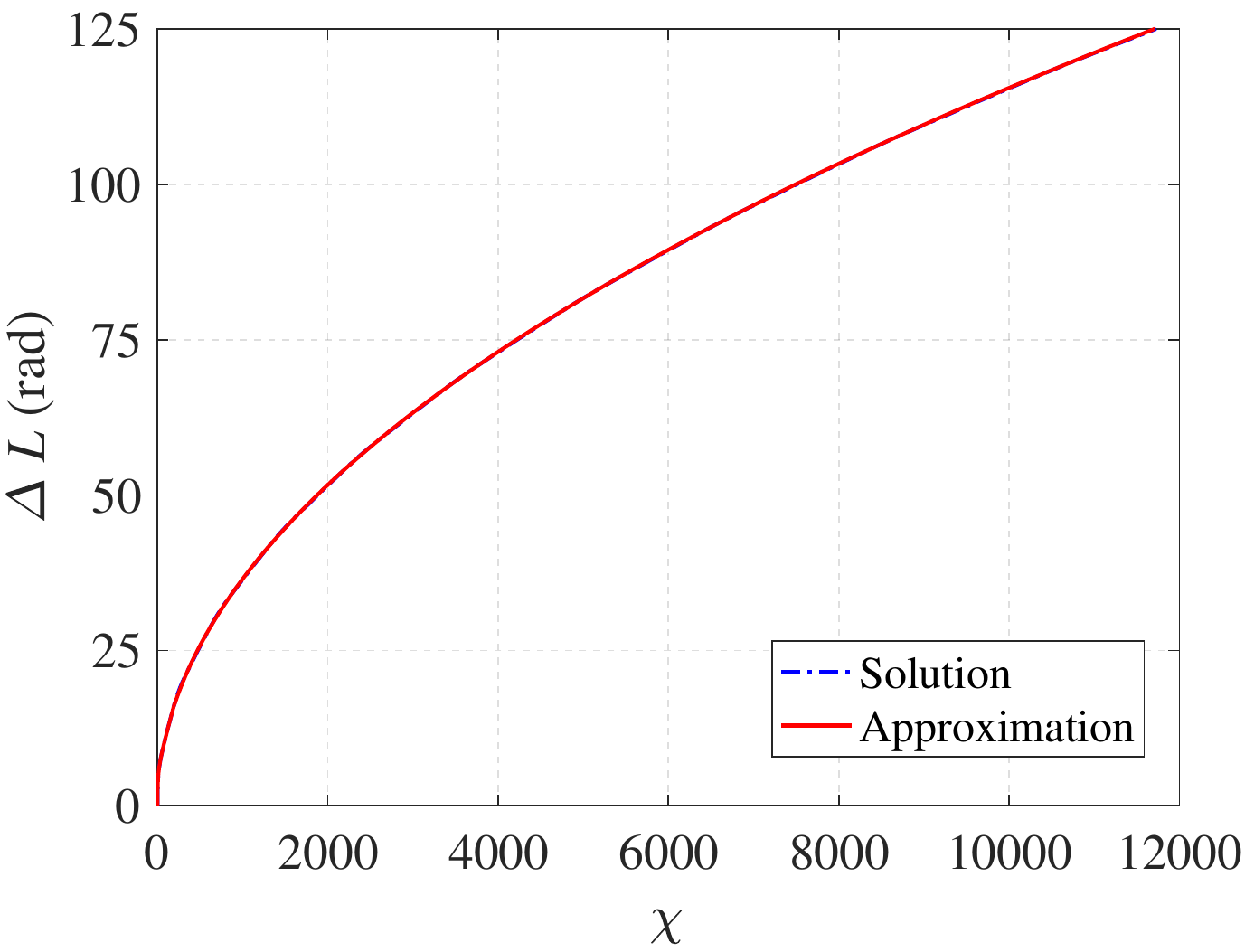}
	\caption{The curves of numerical solution and analytical approximation of $\Delta L$.}
	\label{fig2}
\end{figure}

Based on the solution of $\lambda_1$, the second shooting function can be calculated:
\begin{equation} \label{numerical2}
f_2\left(\Delta L, \, \lambda_1\right) = 2\int_{0}^{\Delta L/2} \frac{ 9\,L^2 + 4 - 6\,\lambda_1\,L\,\sin L - 2\,\lambda_1\,\cos L}{\sqrt{\left(3\, L - 2\,\lambda_1 \sin L\right)^2 + \left(\lambda_1 \cos L - 2\right)^2}} \diff \, L = \chi
\end{equation}
To obtain the analytical Jacobian associated with the shooting function~\eqref{shooting1}, the partial derivatives of $f_2\left(\Delta L, \, \lambda_1\right)$ with respect to $\Delta L$ and $\lambda_1$ are derived as
\begin{equation} \label{partial2}
\begin{aligned}
&\frac{\partial \, f_2}{\partial \, \Delta L} = \frac{9\,\Delta L^2 + 16 - 12\,\lambda_1\,\Delta L\,\sin \left(\Delta L/2\right) - 8\,\lambda_1\,\cos \left(\Delta L/2\right)}{2\,\sqrt{\left[3\,\Delta L - 4\,\lambda_1 \sin \left(\Delta L/2\right)\right]^2 + 4 \left[\lambda_1 \cos \left(\Delta L/2\right) - 2\right]^2}}\\[0.3cm]
&\frac{\partial \, f_2}{\partial \, \lambda_1} = \int_{0}^{\Delta L/2} \frac{ - 2 \, \lambda_1 \, \left(3 \, L \, \cos L - 4\,\sin L\right)^2}{\sqrt{\left(3\, L - 2\,\lambda_1 \sin L\right)^2 + \left(\lambda_1 \cos L - 2\right)^2}^ 3} \diff \, L = 2\,\lambda_1 \frac{\partial \, f_1}{\partial \, \lambda_1}
\end{aligned}
\end{equation}
Combining Eqs.~\eqref{partial1} and~\eqref{partial2}, the analytical Jacobian can be efficiently evaluated by integrating the equation of $\partial f_1 / \partial \lambda_1$. Then, for each pair of $\Delta L$ and $\lambda_1$ displayed in Fig.~\ref{fig1}, the parameter $\chi$ is proved to be positive by calculation. The numerical relationship between $\Delta L$ and $\chi$ is plotted in Fig.~\ref{fig2}. Larger true longitude difference $\Delta L$ is required as the parameter $\chi$ increases. A piecewise function of three segments is developed to fit the numerical solution:
\begin{equation}\label{DeltaL1}
\Delta L \approx \left\{\begin{aligned}
&\; 2 \, \sqrt{\chi} \qquad\qquad\qquad\qquad\quad\; \textnormal{if} \quad \chi \leq 0.2 \\[0.2cm]
&\frac{p_1 \, \chi^3 + p_2 \, \chi^2 + p_3 \, \chi + p_4}{\chi^2 + q_1 \, \chi + q_2} \quad \textnormal{if} \quad 0.2 < \chi \leq 200 \\[0.2cm]
&\; 2 \, \sqrt{\chi \,/\, 3} \qquad\qquad\qquad\quad\;\;\; \textnormal{if} \quad \chi > 200
\end{aligned}\right.
\end{equation}
where $p_1 = 0.04978, \, p_2 = 7.48, \, p_3 = 50.08, \, p_4 = 6.73, \, q_1 = 14.49$, and $q_2 = 15.94$. Figure~\ref{fig2} shows the approximation of the piecewise function~\eqref{DeltaL1} to the numerical solutions with a maximum relative error $\left|\Delta L_{\textnormal{approx}} - \Delta L\right|/\Delta L$ of 0.01. If the true longitude difference is fixed, the maximum final time difference (or rephasing phase) can be approximated by
\begin{equation} \label{chi1}
\left\{
\begin{aligned}
&\chi = \Delta L^2/4, \qquad\qquad\qquad\qquad\qquad\qquad\;\;\, \textnormal{if} \quad \Delta L \leq 0.89 \\[0.3cm]
&\begin{aligned}p_1\, \chi^3 & + \left(p_2 - \Delta L\right)\chi^2 + \left(p_3 - q_1 \Delta L\right)\chi \\[0.2cm]
&+ p_4 - q_2 \Delta L = 0\end{aligned}, \quad\,\, \textnormal{if} \quad 0.89 < \Delta L \leq 16.33 \\[0.3cm]
&\chi = 3 \, \Delta L^2/4, \qquad\qquad\qquad\qquad\qquad\qquad  \textnormal{if} \quad \Delta L > 16.33 \end{aligned}\right.
\end{equation}
Therefore, analytical approximations to the optimal time of flight and the maximum rephasing phase can be obtained by Eqs.~\eqref{DeltaL1} and~\eqref{chi1}, respectively.

%
%

To solve the shooting function $\V{\varPhi}\left(\V{z}\right) = 0$, a convenient strategy is to employ the nonlinear solver Minpack-1 programming package~\cite{more1980user} where the Jacobian is calculated by a forward-difference approximation. The accuracy of computing the Jacobian can be improved by the analytical equations~\eqref{partial1} and~\eqref{partial2}. Another strategy is to use a double loop, in which the inner one solves Eq.~\eqref{numerical1} for $\lambda_1$ and the outer one solves Eq.~\eqref{numerical2} for $\Delta L$, considering the monotone property of these two equations. The analytical approximations~\eqref{lambda_g1} and~\eqref{DeltaL1} serve as initial guesses to the shooting variables. By randomly choosing the parameter $\chi$ in $\left[1.0\times 10^{-5}, \, 1.2 \time 10^{4}\right]$, the accuracy of the initial guesses is then tested by solving 100000 cases with the first strategy and evaluating the number of converged cases and iterations. The solver Minpack-1 is used where the input parameter ``factor'' is set to 0.01, and the differential equations are integrated by the Runge-Kutta adaptive step-size integrator $ode45$ where the relative and absolute tolerances are both set to $10^{-13}$. All the cases converge within 6 iterations on average. The maximum number of iterations required by the solving algorithm is 12. Therefore, the analytical approximations can provide good initial guesses and help the solver converge in several iterations.

\subsection{C. Analytical Solutions to Short- and Long-Term Problems}

When the true longitude difference $\Delta L$ is high-order smaller or larger than 1, some simplifications can be introduced to obtain analytical solutions. In the previous study~\cite{gonzalo2017optimal}, the radial thrust acceleration $a_r$ was identified as zero, while the transversal acceleration $a_{\theta}$ had a bang-bang structure and takes the values $a_{\max}$ or $-a_{\max}$. However, the previous solution obviously violates the boundary constraint $\Delta g = 0$. This work will revisit the short- and long-term problems in order to solve them more accurately.

\textit{1. short-term rephasing}

If the true longitude difference $\Delta L$ is high-order smaller than 1, the function $f_1$ can be transformed into
\begin{equation} \label{short1t}
f_1 \approx \int_{0}^{\Delta L/2} \frac{5\,L^2 - 3\,\lambda_1\,L^2+ 2 - \lambda_1}{\sqrt{\left(3-2\,\lambda_1\right)^2 L^2 + \left(\lambda_1 - 2 - \lambda_1\,L^2/2\right)^2}} \diff \,L = 0 
\end{equation} 
According to Fig.~\ref{fig1}, the costate $\lambda_1$ is about 2, and it is assumed to hold the form $\lambda_1 = 2 - \alpha$, where $\alpha$ is high-order smaller than $\Delta L$. Equation~\eqref{short1t} yields
\begin{equation} \label{short2t}
f_1 \approx \int_{0}^{\Delta L/2} \frac{\alpha - L^2}{\sqrt{\alpha^2 + L^2}}\diff \,L  \approx \alpha \,\textnormal{asinh} \left(\frac{\Delta L}{2\left|\alpha\right|}\right) - \frac{\Delta L^2}{8} = 0
\end{equation}
where $\textnormal{asinh}$ is the arc-hyperbolic sine function. An transcendental equation is then formulated as $\sinh y = 4\,y/\Delta L$, where $y = \Delta L^2 / \left(8\,\alpha\right)$. Based on the numerical solution to this equation, the value of $y$ is larger than 3 when $\Delta L < 1.0$, and the value of $\alpha$ is smaller than $\Delta L^2 / 24$, which confirms the numerical conclusion that $\alpha$ is high-order smaller than $\Delta L$. The exact solution to $\alpha$ (or $\lambda_1$) should be obtained numerically. Nevertheless, the optimal control equation~\eqref{optC} implies that $a_{r}$ is much smaller than $a_{\theta}$ when the true longitude $L$ is not near zero. Then, based on Eq.~\eqref{short2t}, the function $f_2$ can be transformed into
\begin{equation} \label{short3t}
\begin{aligned}
f_2 &\approx 2 \, \int_{0}^{\Delta L/2} \frac{2\,\alpha - L^2}{\sqrt{\alpha^2 + L^2}} \diff \,L \approx 4\,\alpha \textnormal{asinh} \left(\frac{\Delta L}{2\left|\alpha\right|}\right) - \frac{\Delta L^2}{4} \\[0.2cm]
& = \frac{\Delta L^2}{2} - \frac{\Delta L^2}{4} = \chi \end{aligned}
\end{equation}
Therefore, the true longitude difference has the analytical solution $\Delta L = 2\,\sqrt{\chi}$, which holds the same form as the result in Ref.~\cite{gonzalo2017optimal}. By comparison, the optimal control obtained in this work is more accurate and meet all the boundary constraints.

\textit{2. long-term rephasing}

If the true longitude difference $\Delta L$ is high-order larger than 1, the function $f_1$ can only be solved numerically, while the function $f_2$ can be transformed into
\begin{equation} \label{short4t}
f_2 \approx 2 \, \int_{0}^{\Delta L/2} 3\,L \, \diff\,L - \chi = \frac{3\,\Delta L^2}{4} - \chi = 0
\end{equation}
The true longitude difference has the analytical solution $\Delta L = 2\,\sqrt{\chi/3}$, in good agreement with the previous result~\cite{gonzalo2017optimal}.

\section{IV. Propellant-optimal Rephasing Solutions}

Based on the time-optimal solutions, the propellant-optimal rephasing problem will be investigated, where the true longitude difference $\Delta L$ and the parameter $\chi$ are intentionally set to some specific values. The true longitude difference $\Delta L$ should be larger than its minimum value obtained by Eq.~\eqref{DeltaL1} such that the propulsion system can be shut down in part of the trajectory to save propellant. More propellant can be saved with larger $\Delta L$ for a fixed $\chi$. Similar to the time-optimal control solutions, a reduced two-dimensional shooting function will be established based on the first-order optimality conditions. For different values of $\Delta L$ and $\chi$, the contour maps of the two shooting variables and propellant consumption are plotted using numerical solutions. Finally, the analytical propellant-optimal solutions for the cases with relatively small and large true longitude differences $\Delta L$ will be developed.

\subsection{A. Reduced Shooting Function}

The Hamiltonian of the propellant-optimal control problem is established as
\begin{equation} \label{Ham2}
\begin{aligned}
\Ham &= a_{r} \left(\lambda_{\Delta f} \sin L - \lambda_{\Delta g} \cos L\right) + 2 \, a_{\theta} \left(\lambda_{\Delta p} + \lambda_{\Delta f} \cos L + \lambda_{\Delta g} \sin L\right) \\[0.2cm]
&+1.5 \, \lambda_{\Delta t} \, \Delta p - 2\,\lambda_{\Delta t} \, \Delta f\, \cos L - 2\, \lambda_{\Delta t} \, \Delta g \, \sin L+ a \end{aligned}
\end{equation}
which is obtained by replacing the last term of Eq.~\eqref{Ham1} with the thrust acceleration magnitude $a$. Thus, the Euler-Lagrange equations hold the same form as Eq.~\eqref{Eul1}. Based on the symmetry properties introduced in Ref.~\cite{pontani2015symmetry}, the costates take the same value as Eq.~\eqref{Sym1} for the propellant-optimal problem. Bearing in mind Eq.~\eqref{Eul1_Solu}, the optimal control direction is firstly obtained:
\begin{equation} \label{optC2}
\begin{aligned}
\quad a_{r}^{\star}\left(L\right) &= \frac{a \, \left(\lambda_1 \cos L - 2\right)\, \textnormal{sign} \left(\lambda_0\right)}{\sqrt{\left(3\, L - 2\,\lambda_1 \sin L\right)^2 + \left(\lambda_1 \cos L - 2\right)^2}} \\[0.3cm]
a_{\theta}^{\star}\left(L\right) &= \frac{a \, \left(3\,L - 2\,\lambda_1 \sin L\right)\, \textnormal{sign} \left(\lambda_0\right)}{\sqrt{\left(3\, L - 2\,\lambda_1 \sin L\right)^2 + \left(\lambda_1 \cos L - 2\right)^2}}
\end{aligned}
\end{equation}
Substituting Eq.~\eqref{optC2} into Eq.~\eqref{Ham2}, the Hamiltonian is linear in the thrust acceleration magnitude. Therefore, the optimal magnitude is given by
\begin{equation} \label{optC3}
a^{\star} = \left\{\begin{aligned}
&\, 0 \quad\;\;\;\; \textnormal{if} \quad \rho > 0 \\[0.2cm]
& a_{\max} \quad \textnormal{if} \quad \rho < 0\end{aligned} \right.
\end{equation}
where $\rho$ is known as the switching function expressed by
\begin{equation} \label{switch}
\rho = 1 - \left|\lambda_0\right| {\sqrt{\left(3\, L - 2\,\lambda_1 \sin L\right)^2 + \left(\lambda_1 \cos L - 2\right)^2}}
\end{equation}
It is assumed that the switching function takes the value of zero only at some finite isolated points and is neglected in Eq.~\eqref{optC3}. The optimal magnitude is of a bang-bang control structure separated by the switching points where $\rho = 0$. The discontinuous bang-bang control results in numerical difficulties in the accurate integration and Newton-type iteration for solving the optimal control problem~\cite{Jiang2012}. The contour maps in the next subsection will show the effect of the bang-bang control more visually. To ameliorate the above issues, the smoothing technique~\cite{taheri2018generic} is used, through which the optimal magnitude is approximated by
\begin{equation} \label{optC1}
a^{\star} = \frac{a_{\max}}{2} \left[1 + \tanh \left(-\frac{\rho}{\epsilon}\right)\right]
\end{equation}
where $\tanh$ is the hyperbolic tangent function, and the parameter $\epsilon$ determines the smoothness of the control profile. The performance of this approximation to the bang-bang control with different values of $\epsilon$ have been discussed in Ref.~\cite{taheri2018generic} and are omitted here for brevity. In this work, two cases with $\epsilon_1 = 0.1$ and $\epsilon_1 = 0.01$ are tested. The control profile in the first case is smoother, while it is closer to the bang-bang control in the second case. The detail comparison will be presented in the next subsection. The final costates and Hamiltonian are free according to the fixed final states and true longitude difference, respectively. Thus, all the first-order optimality conditions are included in the optimal control equations~\eqref{optC2} and~\eqref{optC1} with the switching function~\eqref{switch}.

Given the values of costates $\lambda_0$ and $\lambda_1$, the final states can be obtained by integrating the ordinary differential equations driven by the optimal control. Similar to the time-optimal control problem, a reduced two-dimensional shooting function can be derived as
\begin{equation} \label{shooting2}
\V{\varPhi} \left(\V{z}\right) = \left[\begin{aligned}
& \;\;\;\, \int_{L_0}^{L_f} \frac{a \left(6\,L\,\sin L + 2\, \cos L - \lambda_1 - 3\,\lambda_1 \, \sin^2 L\right)}{\sqrt{\left(3\, L - 2\,\lambda_1 \sin L\right)^2 + \left(\lambda_1 \cos L - 2\right)^2}} \diff \, L \\[0.2cm]
& \int_{L_0}^{L_f} \frac{ a \left(9\,L^2 + 4 - 6\,\lambda_1\,L\,\sin L - 2\,\lambda_1\,\cos L\right)}{a_{\max} \sqrt{\left(3\, L - 2\,\lambda_1 \sin L\right)^2 + \left(\lambda_1 \cos L - 2\right)^2}} \diff \, L - \chi \end{aligned}\right] = \V{0}
\end{equation}
where $\chi = - \, \textnormal{sign} \left(\lambda_0\right) \Delta t_f \,/\, a_{\max} = \left| \Delta t \right| \,/\, a_{\max}$ is a positive parameter, and $\V{z} = \left[\lambda_0, \, \lambda_1\right]^{\T}$ is a set of the shooting variables. The shooting function is characterized by two main parameters $\Delta L$ and $\chi$. For each specific $\Delta L$ and $\chi$, two branches of solutions can be obtained as $\left\{\lambda_1, \, \lambda_0, \, \Delta t \right\}$ and $\left\{\lambda_1, \, -\lambda_0, \, -\Delta t \right\}$. Without loss of generality, the costate $\lambda_0$ is set to be positive in the rest of this study.

\subsection{B. Numerical Solutions and Linear Interpolation Approximations}

The propellant-optimal solutions with different values of $\Delta L$ and $\chi$ are obtained by numerically solving the shooting function~\eqref{shooting2} to investigate the properties of the solution space. First, a new scaled parameter $\eta \in \left(0,\,1\right)$ is defined to replace the parameter $\chi$:
\begin{equation} \label{eta}
\chi = \left(1 - \eta^2 \right)\, \chi_{\max}
\end{equation}
where $\chi_{\max}$ is the maximum value for a specific $\Delta L$ and can be evaluated by using the time-optimal solution presented in Sec.~\uppercase\expandafter{\romannumeral3}.B. A two-dimensional traversal method is then used to find the relationships between $\left\{\eta, \, \Delta L\right\}$ and $\left\{\lambda_0, \, \lambda_1\right\}$. The true longitude difference $\Delta L$ is traversed between $0.125$ and $125.0$ in a step of $0.125$, and the parameter $\eta$ is traversed between $0.3$ and $0.9$ in a step of $0.001$. The total number of test cases is $1000 \times 601 = 601000$. In simulation, the solver Minpack-1 and integrator $ode45$ are used with the same parameters as in the previous section. Through trial and error, the initial values of costates are randomly guessed in $\lambda_0 \in \left[0, \, 10/\Delta L\right]$ and $\lambda_1 \in \left[-4, \, 4\right]$, and multiple attempts are usually required to solve the two-dimensional shooting function~\eqref{shooting2}. Although the corresponding analytical Jacobian can be similarly derived, it cannot be computed as efficiently as the time-optimal problem, and its derivation is therefore omitted.

\begin{figure}[ht!]
	\centering
	\subfigure[$\epsilon = 0.1$]{\includegraphics[height = 0.26\textheight]{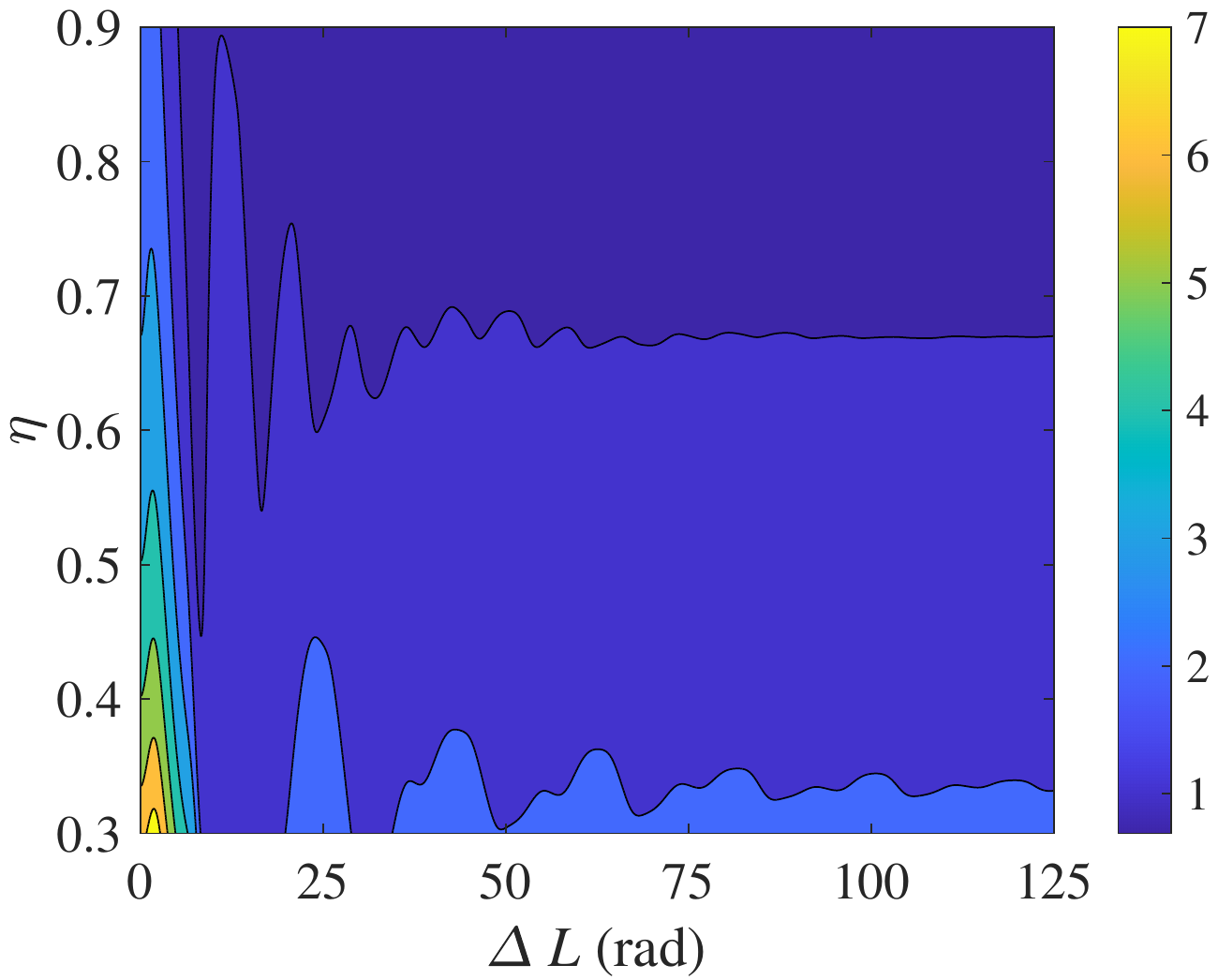}} \;\;
	\subfigure[$\epsilon = 0.01$]{\includegraphics[height = 0.26\textheight]{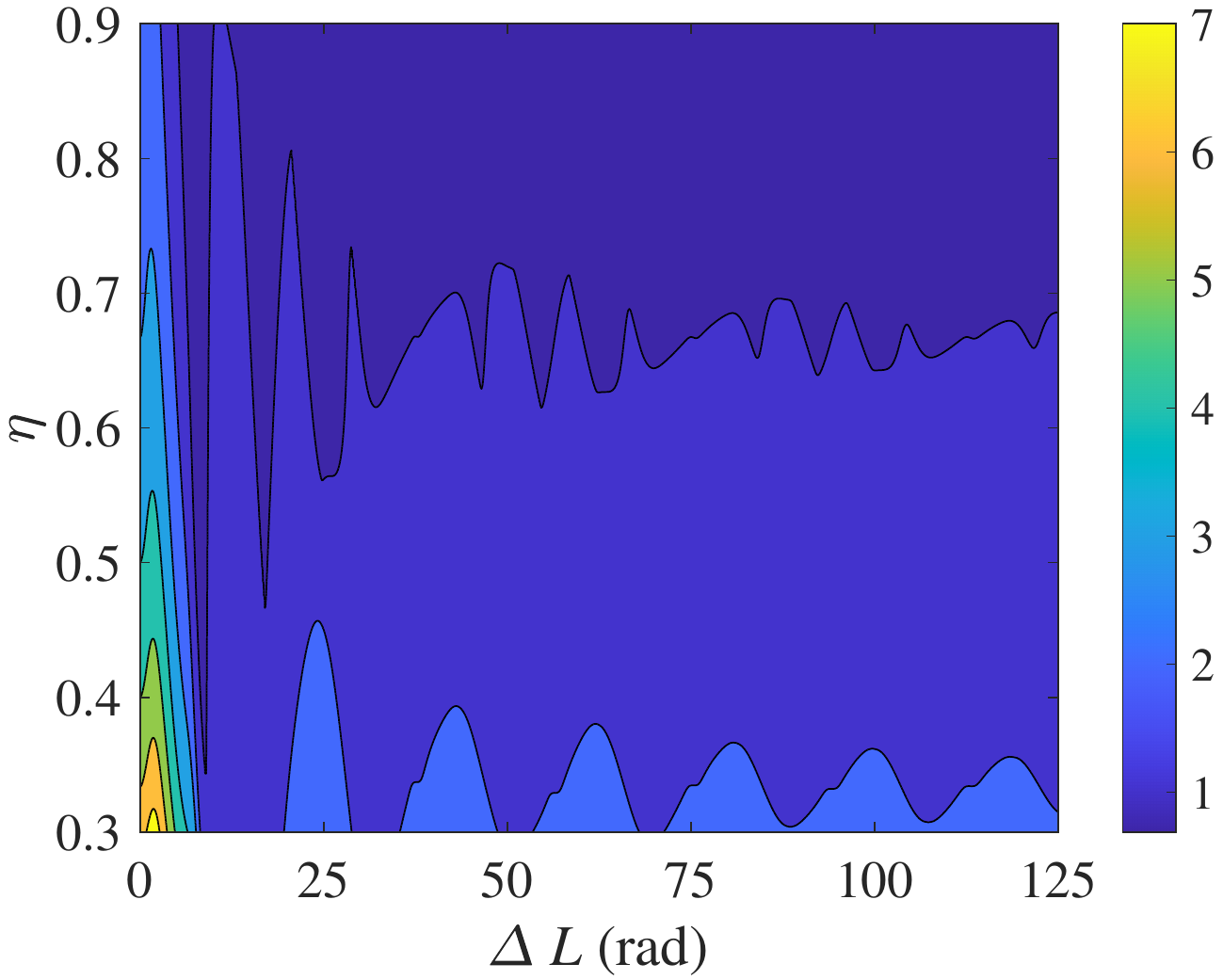}}
	\caption{The results of $\lambda_0 \, \Delta L$ for different values of $\Delta L$ and $\eta$.}
	\label{fig3}
\end{figure}

\begin{figure}[ht!]
	\centering
	\subfigure[$\epsilon = 0.1$]{\includegraphics[height = 0.26\textheight]{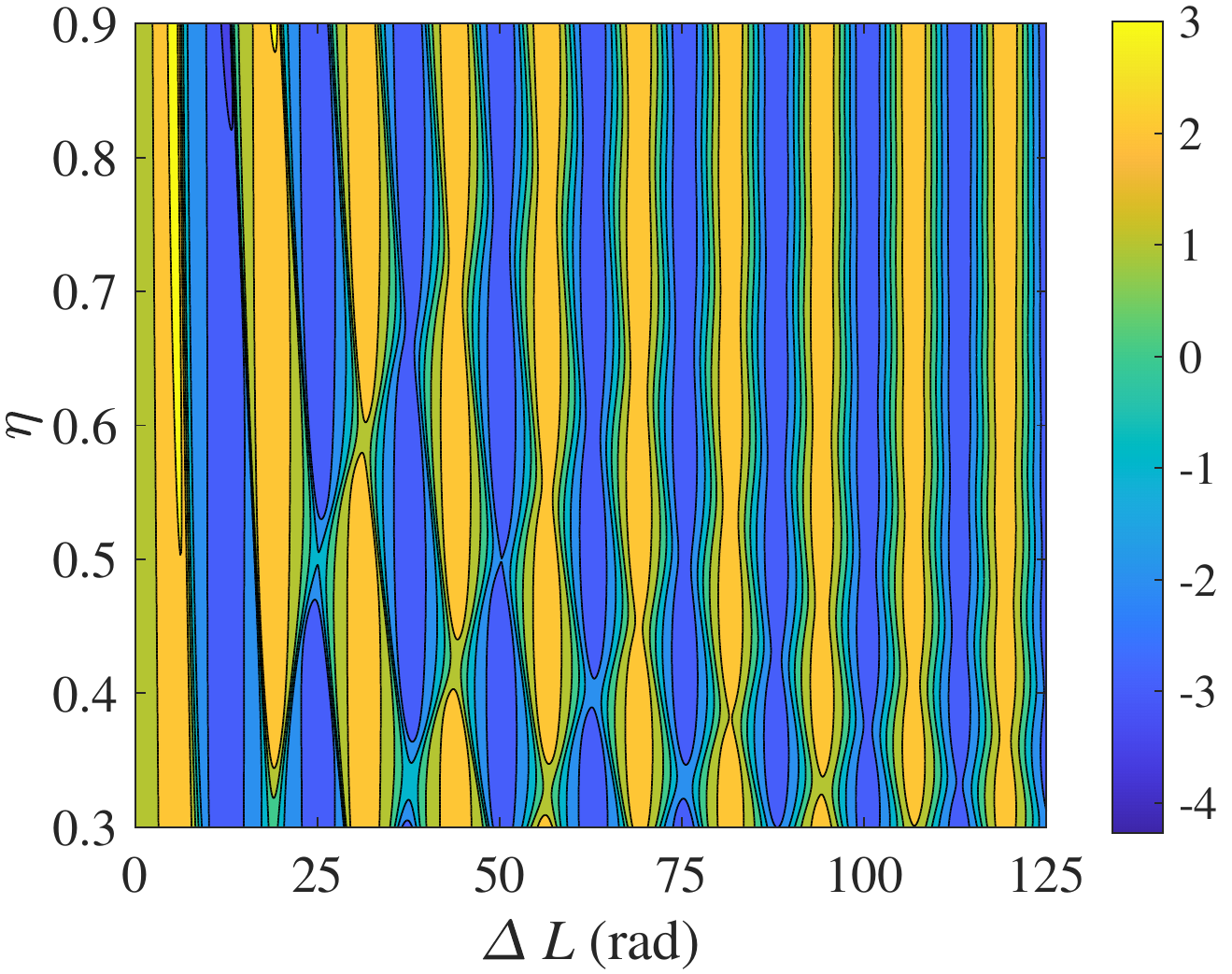}} \;\;
	\subfigure[$\epsilon = 0.01$]{\includegraphics[height = 0.26\textheight]{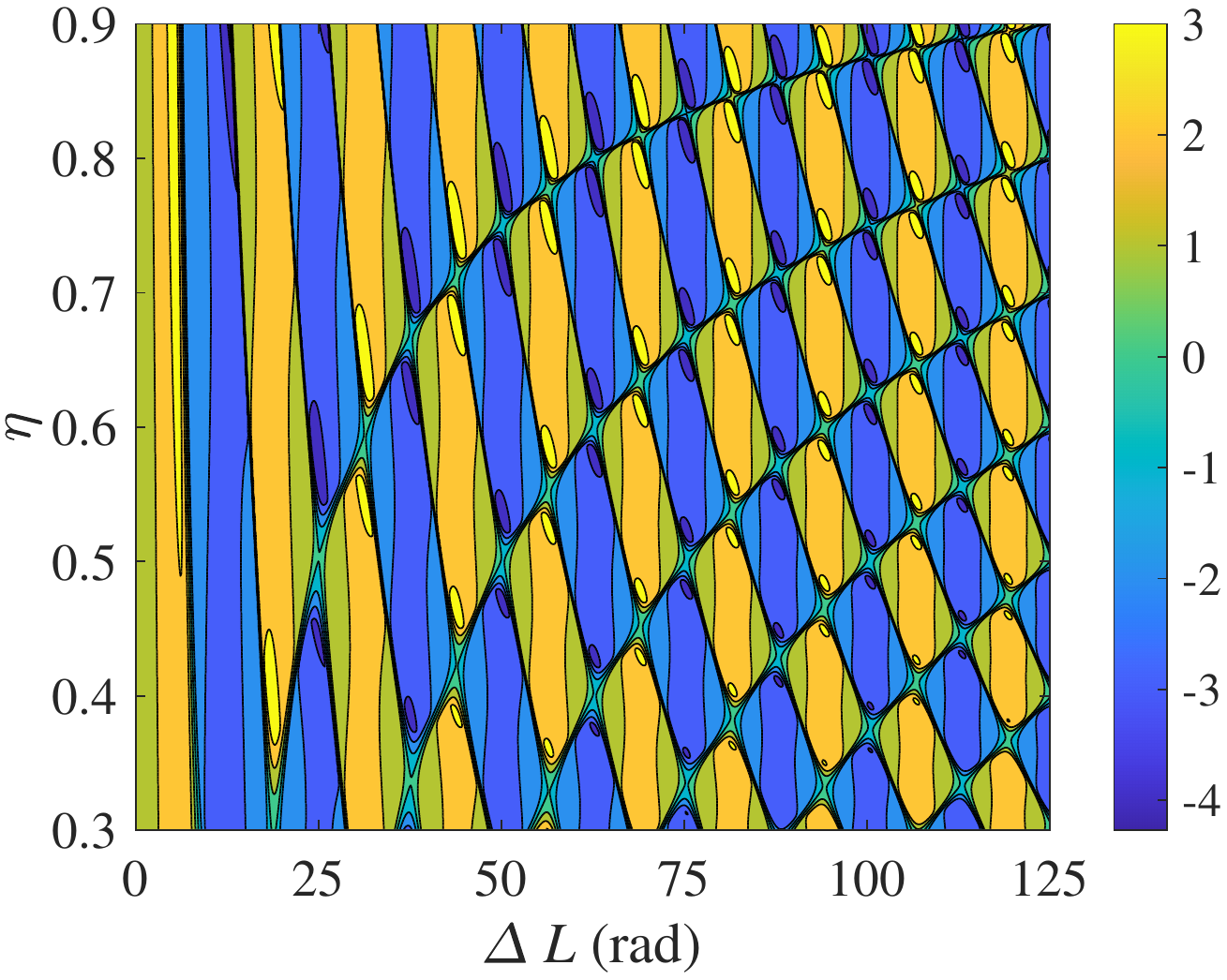}}
	\caption{The results of $\lambda_1$ for different values of $\Delta L$ and $\eta$.}
	\label{fig4}
\end{figure}

\begin{figure}[ht!]
	\centering
	\subfigure[$\epsilon = 0.1$]{\includegraphics[height = 0.258\textheight]{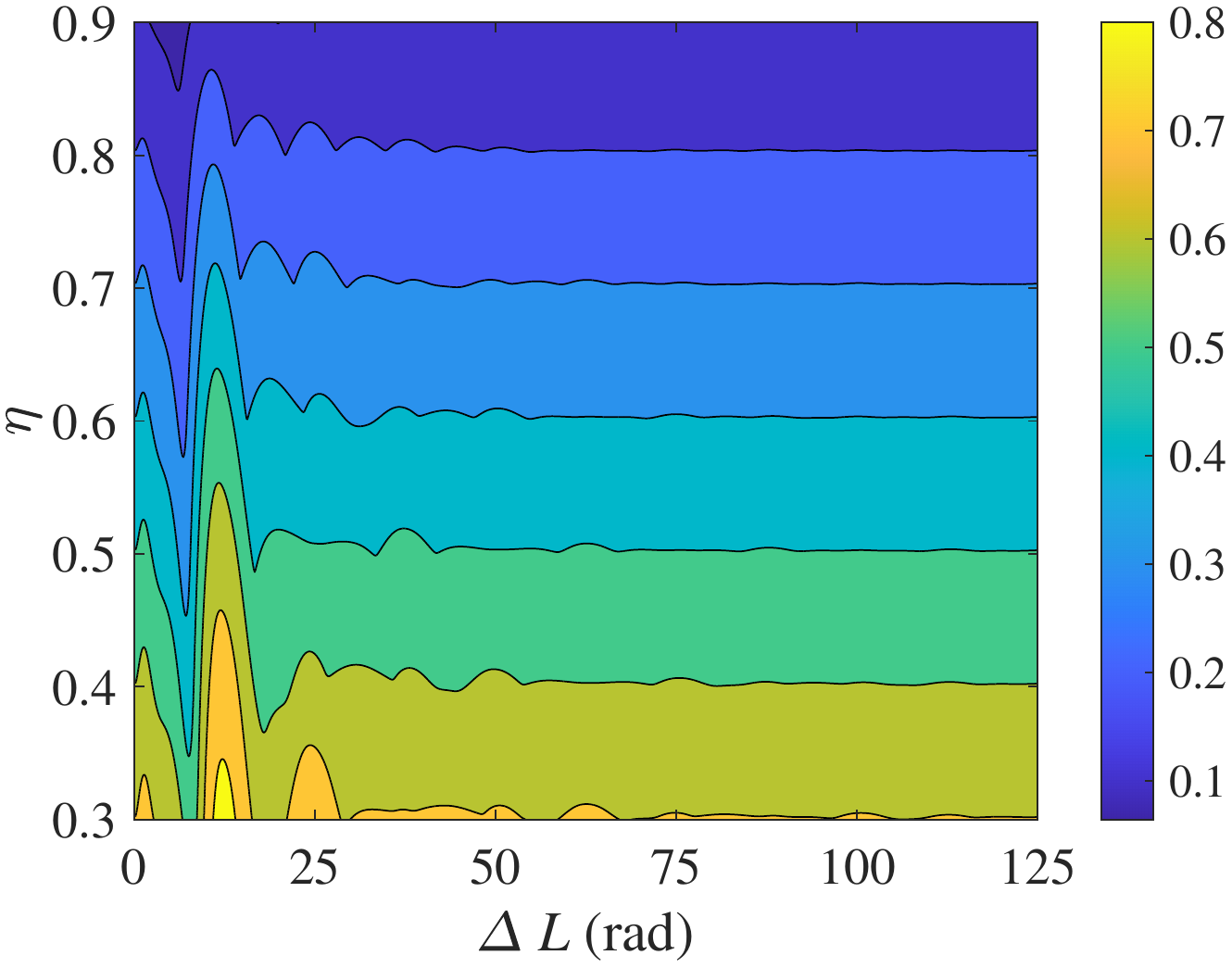}} \;
	\subfigure[$\epsilon = 0.01$]{\includegraphics[height = 0.258\textheight]{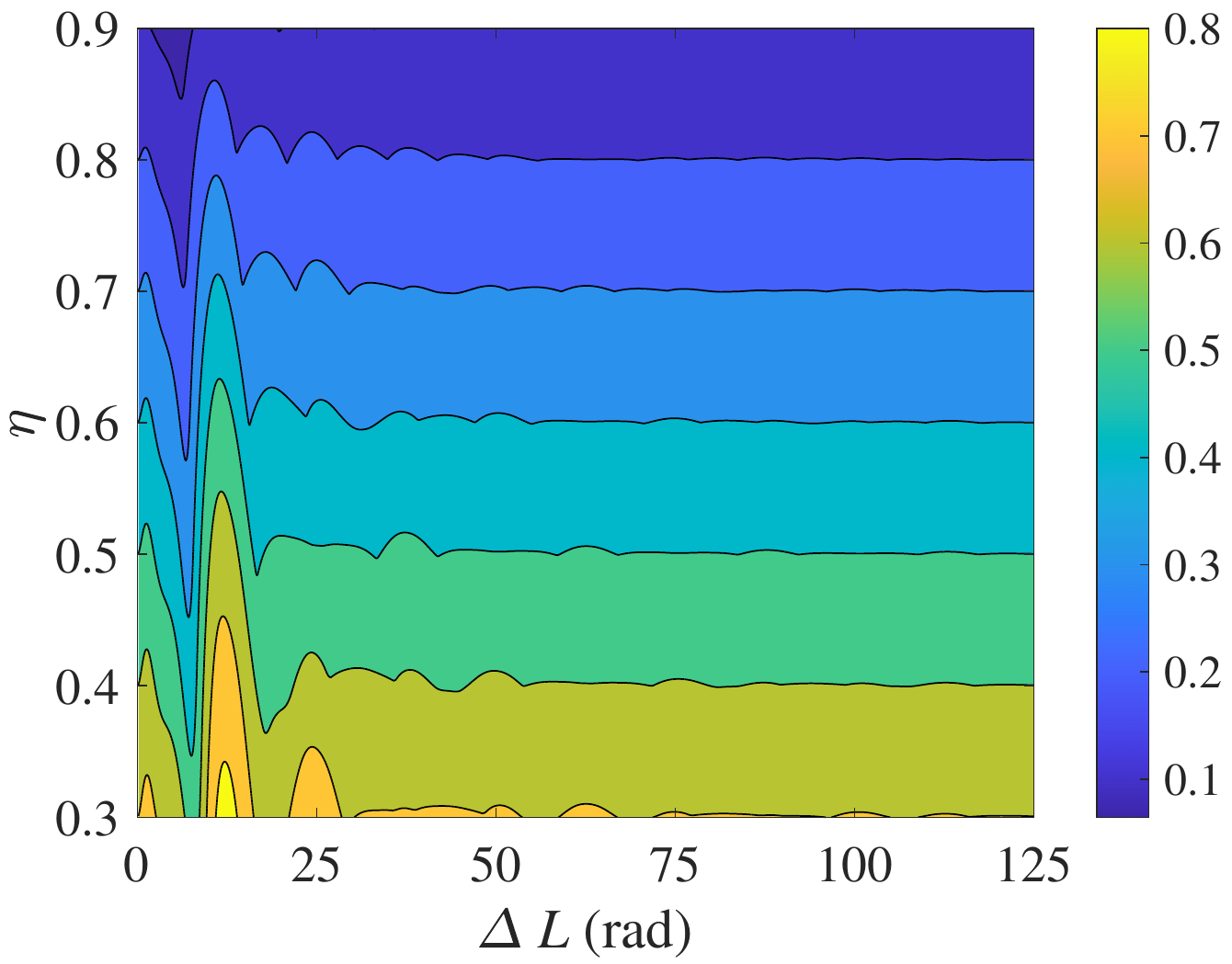}}
	\caption{The results of $J \,/\,\left(a_{\max} \, \Delta L\right)$ for different values of $\Delta L$ and $\eta$.}
	\label{fig5}
\end{figure}

To present the numerical results, the contours of the costates and performance index (or propellant consumption) are given in Figs.~\ref{fig3}--\ref{fig5}. Because $\lambda_0$ decreases rapidly as $\Delta L$ increases, the values of $\lambda_0 \, \Delta L$ rather than $\lambda_0$ are shown in Fig.~\ref{fig3} for clarity. Similarly, the values of $J\,/\,\left(a_{\max} \, \Delta L\right)$ are investigated by Fig.~\ref{fig5}, where $J$ denotes the propellant consumption. Two series of results corresponding to $\epsilon = 0.1$ and $\epsilon = 0.01$ are compared in each figure. As the smoothing parameter $\epsilon$ changes from $0.1$ to $0.01$, the control is closer to the bang-bang control, and the costate values become larger in some cases but smaller in others. It shows that the relationships between the costates $\left\{\lambda_0, \, \lambda_1\right\}$ and the parameters $\left\{\Delta L, \, \eta\right\}$ are notably more complicated when $\epsilon = 0.01$. In some areas of Figs.~\ref{fig3}.b and~\ref{fig4}.b, the costate values vary drastically due to small changes in $\Delta L$ and $\eta$. Figure~\ref{fig5} shows similar propellant consumptions, where the contours on the right are shifted slightly downward and less propellant are therefore consumed. The comparison between the solutions of $\epsilon = 0.01$ and the bang-bang control will be presented later on.

%
%

%
%

%
%

\begin{figure}[ht!]
	\centering
	\includegraphics[scale = 0.56]{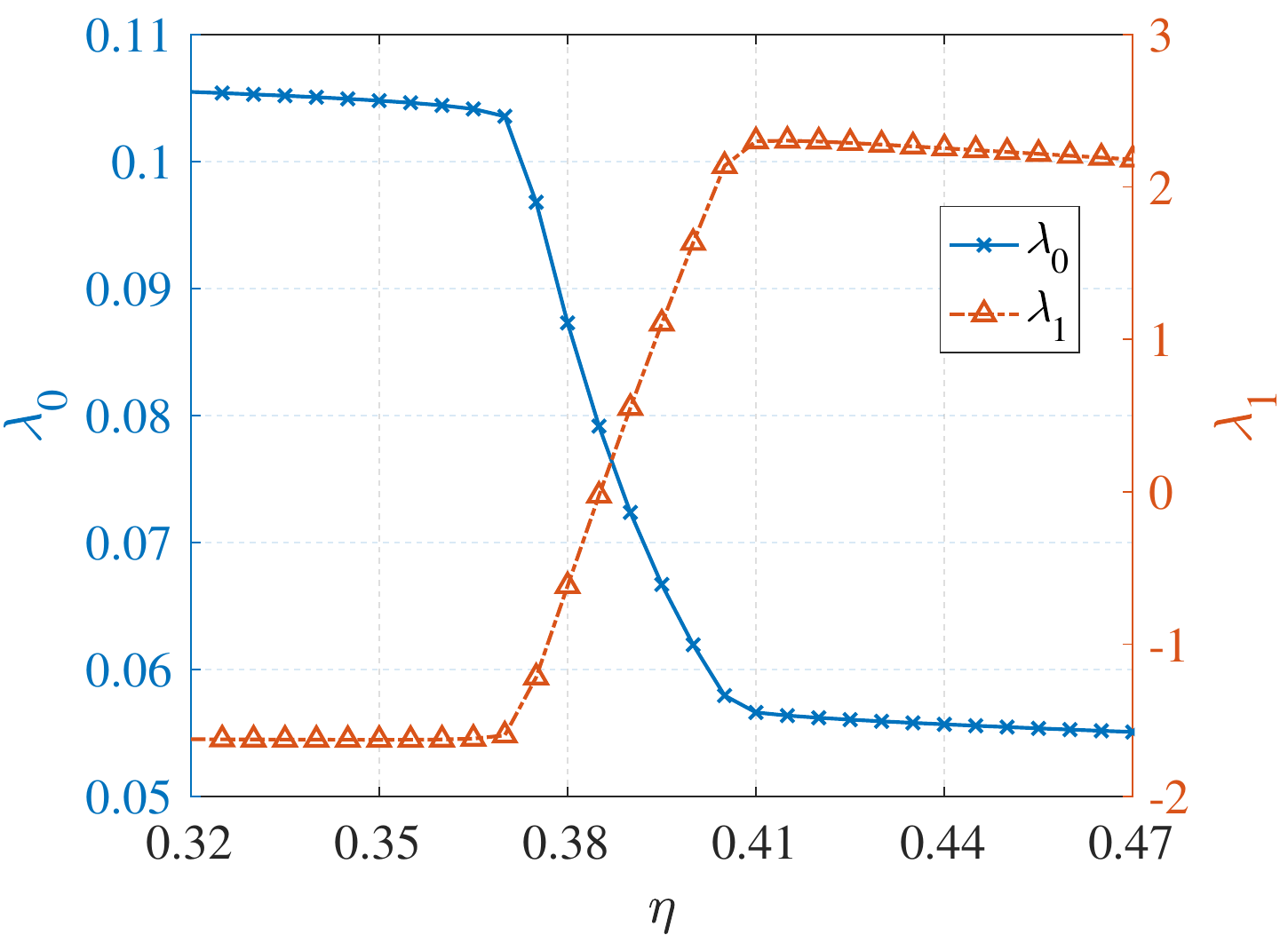} \quad
	\includegraphics[scale = 0.56]{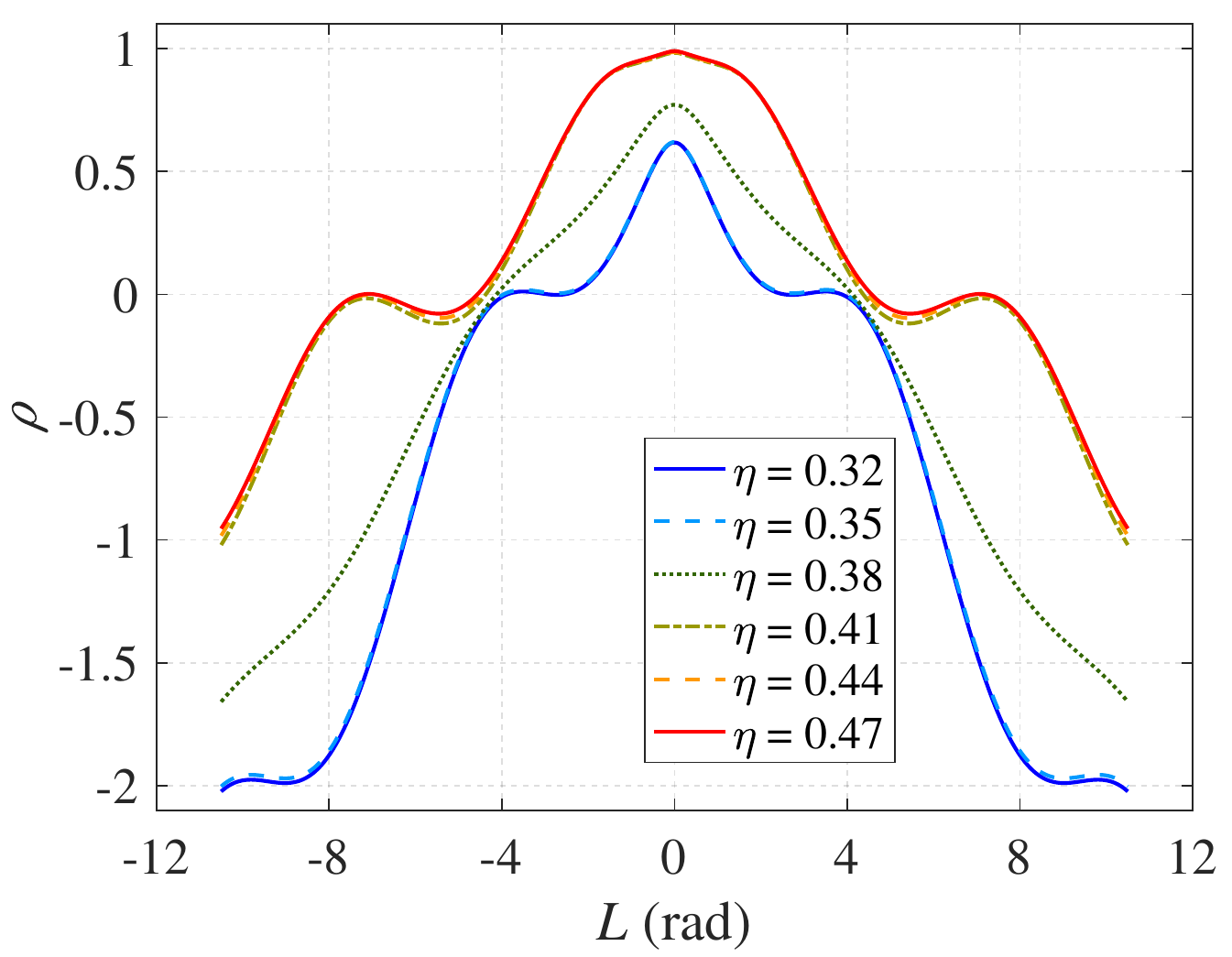}
	\caption{The costates and switching functions as a function of $\eta$ ($\Delta L = 21, \, \epsilon = 0.01$).}
	\label{fig6}
\end{figure}

\begin{figure}[ht!]
	\centering
	\includegraphics[scale = 0.56]{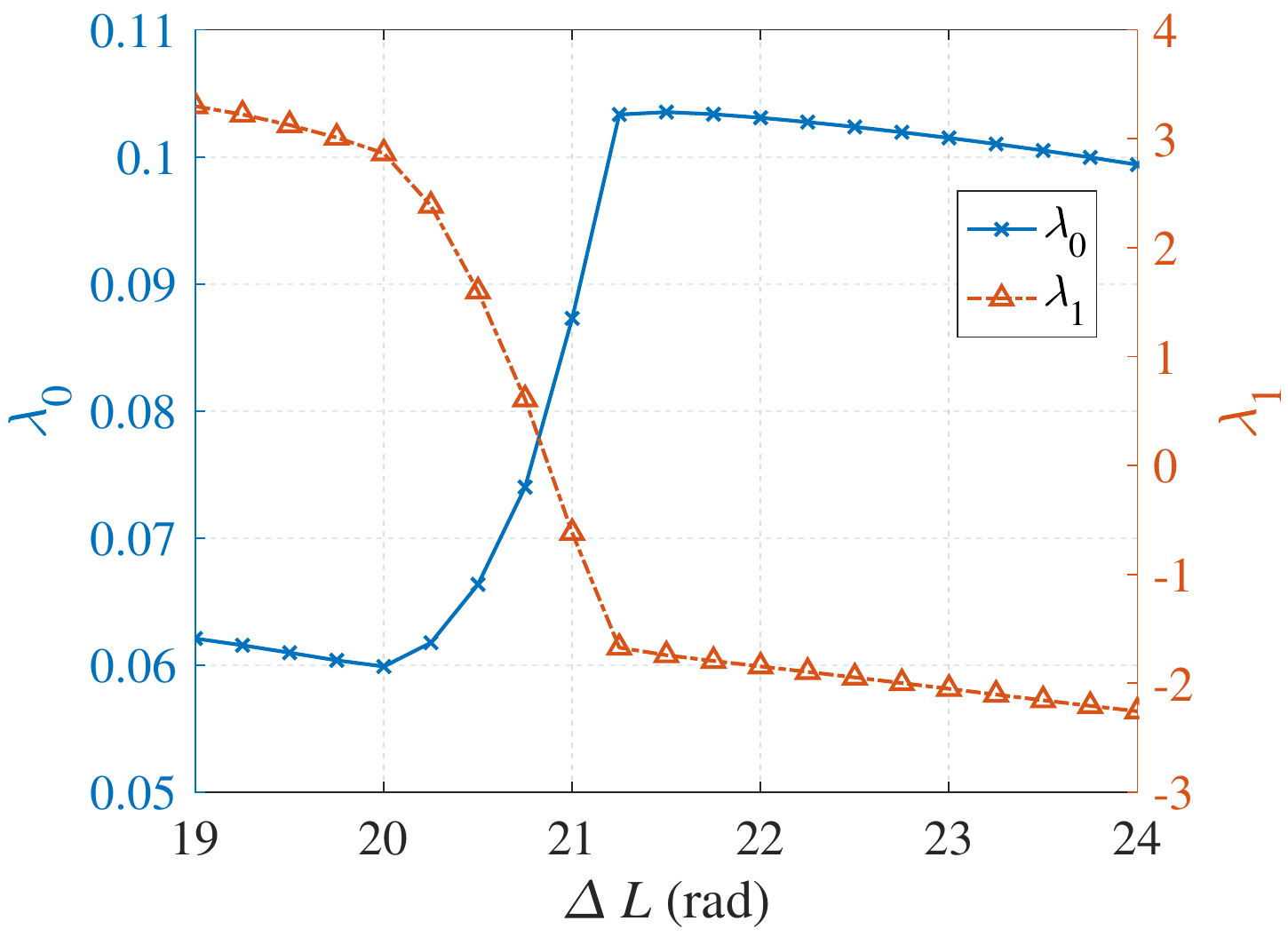} \quad
	\includegraphics[scale = 0.56]{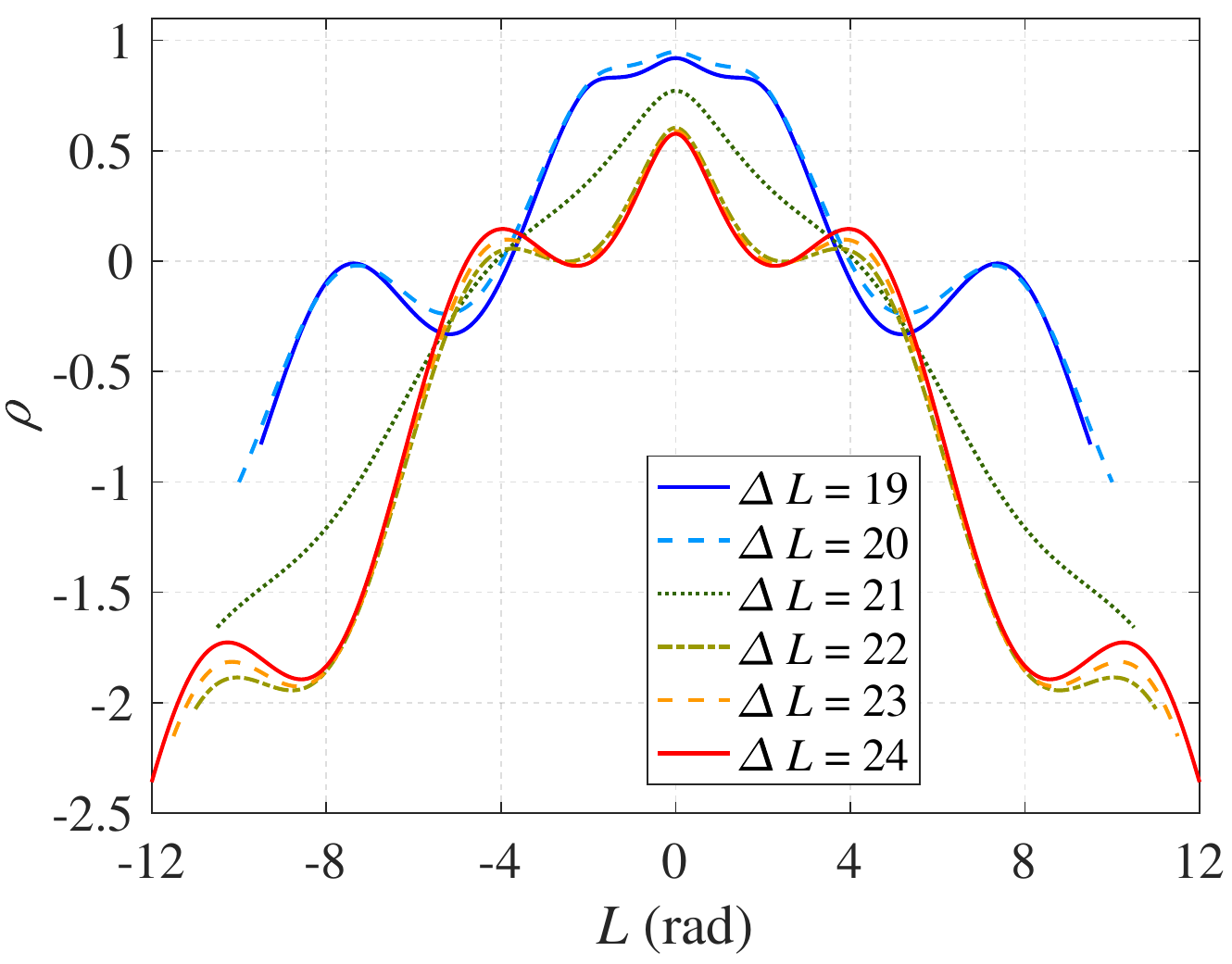}
	\caption{The costates and switching functions as a function of $\Delta L$ ($\eta = 0.38, \, \epsilon = 0.01$).}
	\label{fig7}
\end{figure}

The costates and switching functions $\rho \left(L\right)$ in some cases of Fig~\ref{fig4}.b are investigated to explain the mutations of costate values. For the cases shown in Fig.~\ref{fig6}, the parameter $\eta$ increases from $0.32$ to $0.47$, while the true longitude difference is fixed at $\Delta L = 21$. The mutations occur approximately at $\eta = 0.38$. Taking this result as the demarcation, the switching functions are clearly divided into two clusters. In each cluster, the optimal solutions with $\eta = 0.32$ and $\eta = 0.47$ have four burning arcs, while the others only have two burning arcs. Although the total propellant consumption decreases monotonically according to Fig.~\ref{fig5}.b, the number of burning arcs changes from 4 to 2 and finally to 4 as $\eta$ increases from $0.32$ to $0.47$. Thus, when $\Delta L$ is fixed and $\eta$ varies, the change in the number of burning arcs results in the clustering of switching functions and the mutations of costates. In Fig.~\ref{fig7}, the cases with fixed $\eta$ and increasing $\Delta L$ are presented. Similarly, the switching functions are divided into two clusters, and the mutations occur approximately at $\Delta L = 21$. When both parameters $\Delta L$ and $\eta$ vary, each small positive or negative region in Fig.~\ref{fig4}.b represents a cluster of switching functions. The number of burning arcs might change in each cluster (i.e., 2 or 4 for all numerical solutions), and the costate mutates as the switching function changes from one cluster to another.

%
%

%
%

To efficiently approximate the propellant consumption and costate values, the linear interpretation is used based on the numerical solutions. The approximate costate values serve as the initial guesses to solve the two-dimensional shooting function~\eqref{shooting2}. By randomly choosing the parameter $\Delta L$ in $\left[0.125, \, 125\right]$ and $\eta$ in $\left[0.3, \, 0.9\right]$, the performance of the approximation is then tested by solving 100000 cases and evaluating the number of converged cases and iterations. When the smoothing parameter $\epsilon$ is set to 0.1, all cases converge within 5 iterations on average, and the maximum number of iterations is 14. After setting $\epsilon = 0.01$ and using the corresponding approximation, all cases can converge within 5 iterations on average, but the maximum number of iterations becomes 28. 


\subsection{C. Analytical Solutions to Short- and Long-Term Problems}

For the propellant-optimal control problem, two analytical solutions can be developed when the true longitude difference $\Delta L$ is relatively small and large, respectively. The simplification techniques for these two solutions are introduced as follows.

\textit{1. short-term rephasing}

If the true longitude difference $\Delta L$ is high-order smaller than 1, the costate $\lambda_1$ is about 2 according to the numerical solutions. The expression $\lambda_1 = 2 - \alpha$ is employed, where $\alpha$ is high-order smaller than $\Delta L$. The switching function~\eqref{switch} can be transformed into
\begin{equation} \label{shortm1}
\rho = 1 - \left|\lambda_0\right| \sqrt{\alpha^2 + L^2}
\end{equation} 
which is a monotonic increasing function with respect to $\left| L \right|$. Therefore, the switching function is assumed to be positive at the ranges $\left[-\Delta L/2, -\beta\,\Delta L/2\right]$ and $\left[\beta\,\Delta L/2, \Delta L/2\right]$, where $\beta$ is a positive parameter smaller than 1.

Based on the derivations of Eqs.~\eqref{short1t}--\eqref{short3t}, an analytical solution can be derived as
\begin{equation} \label{shortm2}
\frac{\left(1-\beta^2\right) \Delta L^2}{4} - \chi = 0
\end{equation}
Since $\chi = \left(1-\eta^2\right) \chi_{\max} = \left(1-\eta^2\right) \Delta L^2 / 4$, the parameter $\beta$ equals to the parameter $\eta$, i.e., $\beta = \eta$. Then, the costate $\lambda_0$ is obtained as $\left|\lambda_0\right| \approx 1 \,/\left(\eta \, \Delta L\right)$, and the propellant consumption is computed by
\begin{equation} \label{shortm3}
J = \left(1-\beta\right) a_{\max} \, \Delta L = \left(1-\eta\right) a_{\max} \,\Delta L
\end{equation}
The parameter $J \,/\, \left(a_{\max} \, \Delta L\right) = 1 - \eta$ is linear in $\eta$, which has been confirmed by the numerical results in Fig~\ref{fig5}.

\textit{2. long-term rephasing}

If the true longitude difference $\Delta L$ is high-order larger than 1, the switching function $\rho$ can be estimated by
\begin{equation} \label{short3}
\rho = 1 - 3\left|\lambda_0 \, L \right|
\end{equation}
Similarly, the switching function can be assumed to be positive at the ranges $\left[-\Delta L/2, -\beta\,\Delta L/2\right]$ and $\left[\beta\,\Delta L/2, \Delta L/2\right]$. Based on the derivations of Eq.~\eqref{short4t}, the second shooting function is transformed into
\begin{equation} \label{short2}
\frac{3 \, \left(1-\beta^2\right) \Delta L^2}{4} - \chi = 0
\end{equation}
where the parameter $\chi = \left(1-\eta^2\right) \chi_{\max} = 3\left(1-\eta^2\right) \Delta L^2 / 4$. The similar results of $\beta = \eta$ and $\left|\lambda_0\right| \approx 1 \,/\left(3\,\eta \, \Delta L\right)$ can be obtained. Finally, the propellant consumption is given by
\begin{equation} \label{shortm4}
J = \left(1-\beta\right) a_{\max} \, \Delta L = \left(1-\eta\right) a_{\max} \,  \Delta L
\end{equation}
The values of $J \,/\, \left(a_{\max} \, \Delta L\right)$ are equal in the short- and long-term rephasing problems with the same $\eta$, and this characteristic has been presented in Fig.~\ref{fig5}.

\section{V. Numerical Examples}

In this section, the numerical results of the costate values, optimal control profiles, and trajectories obtained with the linear dynamics~\eqref{dyn4} will be presented and compared with those obtained with the nonlinear dynamics~\eqref{dyn}. The optimal control problem with the nonlinear dynamics is first introduced, and two four-dimensional shooting functions are established for the time- and propellant-optimal problems. Then, several cases are tested to compare the nonlinear and linear results. 

The two nonlinear equations~\eqref{dyn} and~\eqref{dyn2} are equivalent and lead to the same optimal solution. The costates concerned with the indirect methods can be converted to each other according to Ref.~\cite{wu2021minimum}. In this work, the nonlinear equation~\eqref{dyn2} is employed. Bearing in mind the performance index~\eqref{per}, the associated Hamiltonian can be derived as
\begin{equation} \label{nonHam}
\Ham = \left[\V{\lambda}_{x}^{\T} \, \V{B}\left(\V{x}, \, L\right) \V{a} + \lambda_t + \phi \right]\frac{1}{A\left(\V{x},\,L\right)}
\end{equation}
where $\left[\V{\lambda}_x^{\T}, \, \lambda_t\right]^{\T}$ is the costate vector associated with the state $\left[\V{x}^{\T},\,t\right]^{\T}$. The Euler-Lagrange equations are obtained as
\begin{equation} \label{nonEul}
\left\{
\begin{aligned}
\V{\lambda}_x^{\prime} &= -\frac{\partial \Ham}{\partial \V{x}} = -\V{\lambda}_x^{\T} \, \frac{\partial}{\partial \V{x}} \left(\frac{\V{B} \,\V{a}}{A}\right) + \frac{\lambda_t+\phi}{A^2}\frac{\partial A}{\partial \V{x}}\\[0.3cm]
\V{\lambda}_t^{\prime} &= -\frac{\partial \Ham}{\partial t} = 0
\end{aligned}\right.
\end{equation}
where $\phi$ equals to 1 or $a$. Therefore, the costate $\lambda_t$ is constant. Considering that the thrust acceleration is high-order small than ${\mu \,/\, \hat{p}^2} = 1$, the nonlinear Euler-Lagrange equations~\eqref{nonEul} can be transformed into the linear equations~\eqref{Eul1} by setting $\left[\V{\lambda}_x^{\T}, \, \lambda_t + 1\right]^{\T} = \left[\V{\lambda}_{\Delta x}^{\T}, \, \lambda_0 \right]^{\T}$ for the time-optimal problem and $\left[\V{\lambda}_x^{\T}, \, \lambda_t\right]^{\T} = \left[\V{\lambda}_{\Delta x}^{\T}, \, \lambda_0 \right]^{\T}$ for the propellant-optimal one. To minimize the Hamiltonian, the optimal control should follow
\begin{equation} \label{nonoptC}
\V{a}^{\star} = -a^{\star} \frac{\V{B}\,\V{\lambda}_x}{\left\|\V{B}\,\V{\lambda}_x\right\|}
\end{equation}
where $a^{\star}$ takes the value 1 for the time-optimal problem and is of bang-bang control for the propellant-optimal one. By introducing the smoothing technique, the optimal magnitude holds the same form as Eq.~\eqref{optC1}, where the switching function is given by
\begin{equation} \label{nonswitch}
\rho = 1 - \left\|\V{B} \,\V{\lambda}_x\right\|
\end{equation}
These optimal control equations can be transformed into the equations~\eqref{optC} and~\eqref{optC2}--\eqref{optC1} by the linearization technique. Therefore, the proposed solutions can be used to approximate the nonlinear solutions. Based on the analyses in Sec.~\uppercase\expandafter{\romannumeral3}.A, the transversality condition is automatically guaranteed, and the costates at the final true longitude are free. The remaining boundary constraints are written as
\begin{equation} \label{nonboundary}
\V{x}\left(L_f\right) = \CC{x}, \quad t\left(L_f\right) = t_f
\end{equation}
where $t_f$ is the final time when the active satellite rendezvous with the target. Note that the quantity $\left|\lambda_t + \phi\right|$ is scaled to 1 when solving the time-optimal problem.

The two-point boundary value problem is then formulated and solved by a shooting method. The shooting function is expressed by Eq.~\eqref{nonboundary}, and the shooting variables are identified as $\V{z} = \left[\V{\lambda}_x^{\T}, \, \Delta L\right]^{\T}$ and $\V{z} = \left[\V{\lambda}_x^{\T}, \, \lambda_t\right]^{\T}$ for the time- and propellant-optimal control problems, respectively. Once the shooting variables are obtained, the state and costate values at the final time can be evaluated by integrating the nonlinear equations~\eqref{dyn2} and Euler-Lagrange equation~\eqref{nonEul} driven by the optimal control Eq.~\eqref{nonoptC}.

The proposed time-optimal solution in Sec.~\uppercase\expandafter{\romannumeral3} depends on only one key parameter $\chi = -\textnormal{sign}\left(\lambda_0\right) \Delta t_f \,/\, a_{\max}$. By comparison, the solution is affected by both the two parameters $\Delta t$ and $a_{\max}$ in the nonlinear formulation. The final time difference $\Delta t_f$ equals to the phase difference $\Delta \theta$ to be rephased, and $a_{\max}$ is evaluated by the ratio of the thrust acceleration to the gravitational acceleration. For the geocentric rephasing problem, the gravitational acceleration is about $8.37 \, \unit{m/s^2}$ in a low Earth circular orbit, and it is about $0.22 \, \unit{m/s^2}$ in the geostationary orbit. Thus, it is reasonable to assume that $a_{\max}$ is high-order smaller than 1 for the satellites equipped with a low-thrust propulsion system.



\subsection{A. Linear and Nonlinear Time-Optimal Rephasing Solutions}

To illustrate a systematic comparison between the linear and nonlinear time-optimal solutions, three cases with $\chi = 0.05, \, 10, $ and $1000$ are tested, in which the parameters are summarized in Table~\ref{table1}. A wide range of values are taken for the key parameter $\chi$ and thrust acceleration $a_{\max}$, while the final time (or phase) difference $\Delta t_f$ are then computed. The maximum values of the thrust acceleration and phase difference are $0.1 \, \mu \,/\, \hat{p}^2$ and $1\, \unit{rad}$, respectively. In Case 2, three nonlinear solutions characterized by the same key parameter and different thrust acceleration are compared.

%
{
\begin{table}[!htb]
	\caption{Parameters of the linear and nonlinear time-optimal solutions}
	\label{table1}
	\vspace{0.25cm}
	\centering \small
	\begin{tabular}{cccccc}
		\hline \hline
		\multirow{2}{*}{Parameter} & \multirow{2}{*}{Case 1} & \multicolumn{3}{c}{Case 2} & \multirow{2}{*}{Case 3} \\ 
		 & & S1 & S2 & S3 & \\ \hline
		$\chi$     & 0.05  & \multicolumn{3}{c}{10}  & 1000 \\
		$\Delta t_f$ & $-0.005$ & $-0.01$   & $-0.1$  & $-1.0$ & $-1.0$ \\
		$a_{\max}$ & 0.1   & 0.001 & 0.01 & 0.1 & 0.001 \\
		\hline \hline
	\end{tabular}
\end{table}}
%

%
{
\begin{table}[!htb]
	\caption{Comparison of the linear and nonlinear time-optimal solutions}
	\label{table2}
	\vspace{0.25cm}
	\centering \small
	\begin{tabular}{cccc}
		\hline \hline
		Case & Solution & $\V{\lambda}_x$ & $\Delta L$ \\ \hline
		\multirow{2}{*}{1} & Linear & $\left[0.33650,\, -0.44491,\, 0.04464\right]^{\T}$ & 0.44866 \\
		                & Nonlinear & $\left[0.33160,\, -0.43755,\, 0.04477\right]^{\T}$ & 0.45366 \\
		\multirow{4}{*}{2} & Linear & $\left[3.75470,\, -1.19191,\, 3.70636\right]^{\T}$ & 5.00627 \\
		                       & S1 & $\left[3.74128,\, -1.17964,\, 3.69340\right]^{\T}$ & 5.01167 \\
		                       & S2 & $\left[3.62345,\, -1.07286,\, 3.57886\right]^{\T}$ & 5.06025 \\
		                       & S3 & $\left[2.68055,\, -0.28892,\, 2.61255\right]^{\T}$ & 5.55308 \\
		\multirow{2}{*}{3} & Linear & $\left[27.30648,\, 1.20278,\, -1.06349\right]^{\T}$ & 36.40864 \\
		                & Nonlinear & $\left[26.18922,\, 0.48488,\, -2.17276\right]^{\T}$ & 37.19677 \\
		\hline \hline
	\end{tabular}
\end{table}}
%

Firstly, the linear solution is obtained by the method presented in Sec.~\uppercase\expandafter{\romannumeral3}, and the costate $\V{\lambda}_{\Delta x}$ is computed by Eq.~\eqref{Eul1_Solu}. Then, the nonlinear solution is obtained by solving the four-dimensional shooting function~\eqref{nonboundary} with the initial guess given by the linear solution. All cases converge within several iterations, and the results are presented in Table~\ref{table2}. Comparing the results of Case 1 with S3 (or Case 3 with S1), the linear solutions are closer to the nonlinear ones when the thrust acceleration is the same and the phase difference is smaller. Comparing the results of Case 3 with S3, the linear solutions are closer to the nonlinear ones when the thrust acceleration is smaller and the phase difference is the same. The relative error of $\Delta L$ takes its maximum value of $9.8\%$ in S3 and $0.1\%$ in S1, and the relative errors in the other solutions are about $1.5\%$. Therefore, the linear solutions provide good initial guess to the nonlinear solutions and accurate approximation to the optimal time of flight for the low-thrust rephasing problem. The piecewise function~\eqref{DeltaL1} is suggested to analytically approximate the optimal time of flight.

%
\begin{figure}[ht!]
	\centering
	\includegraphics[scale = 0.56]{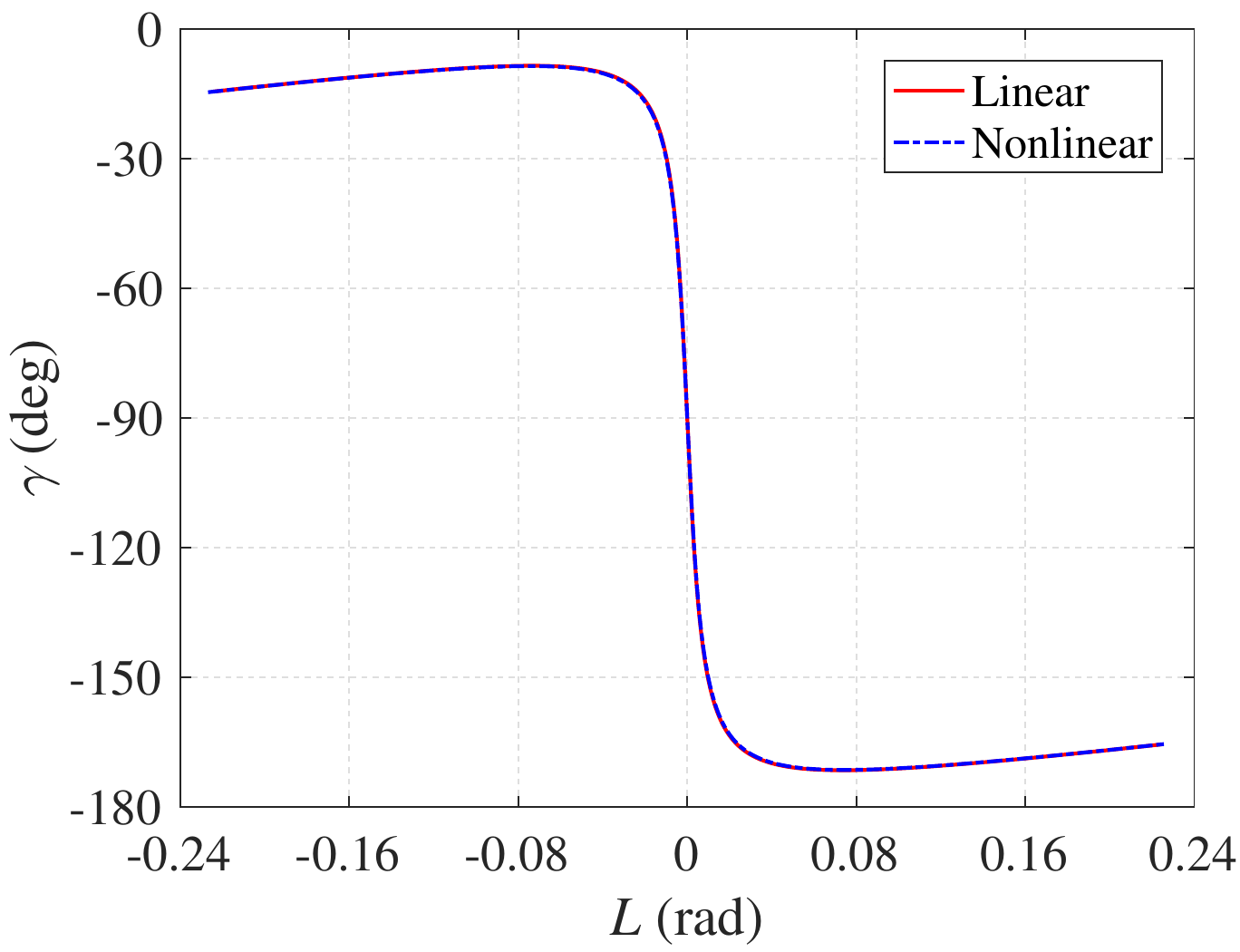} \quad
	\includegraphics[scale = 0.57]{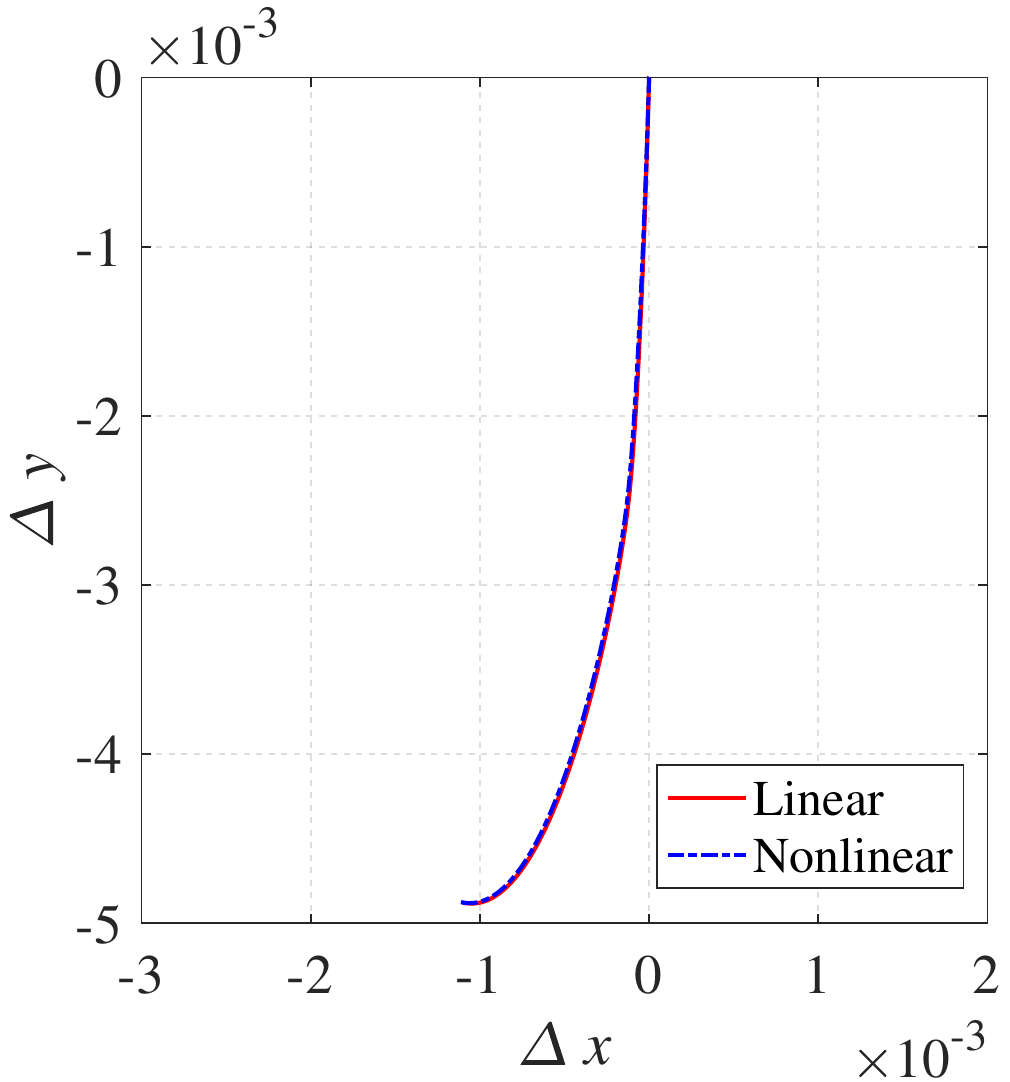} \qquad
	\caption{The control profiles and relative trajectories in $xy$ plane in Case 1.}
	\label{fig8}
\end{figure}
%

%
\begin{figure}[ht!]
	\centering
	\includegraphics[scale = 0.56]{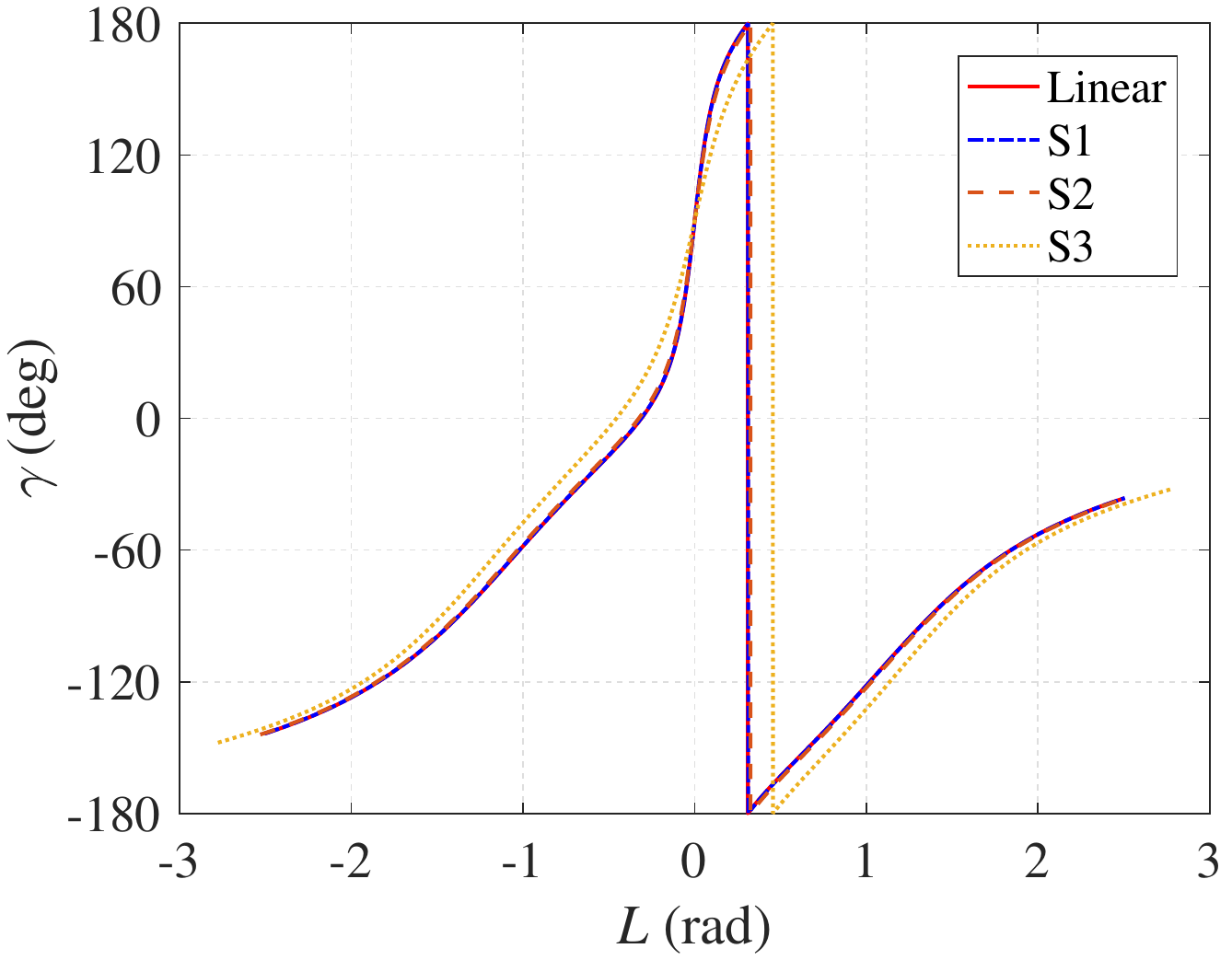} \quad
	\includegraphics[scale = 0.56]{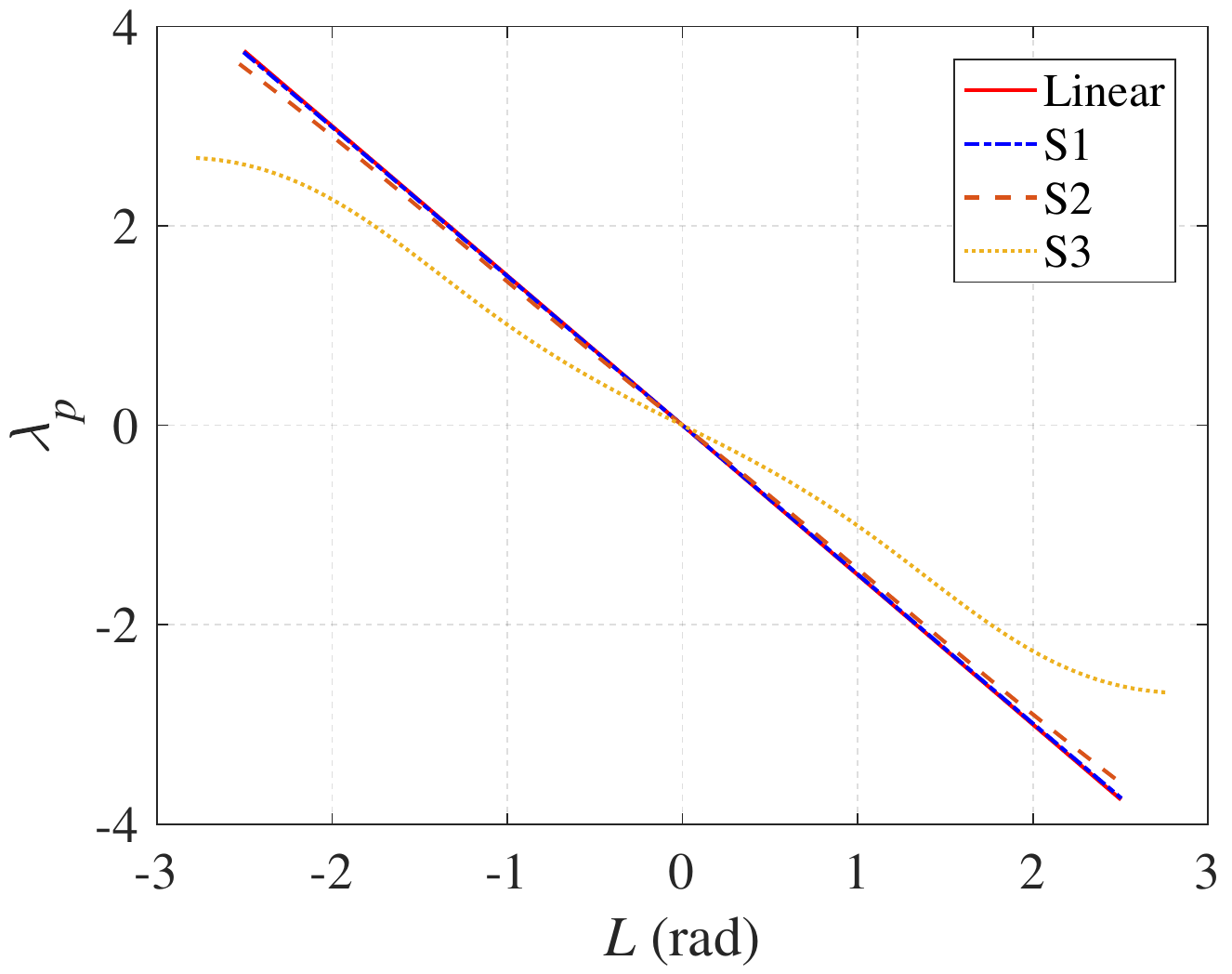} \\
	\;\includegraphics[scale = 0.56]{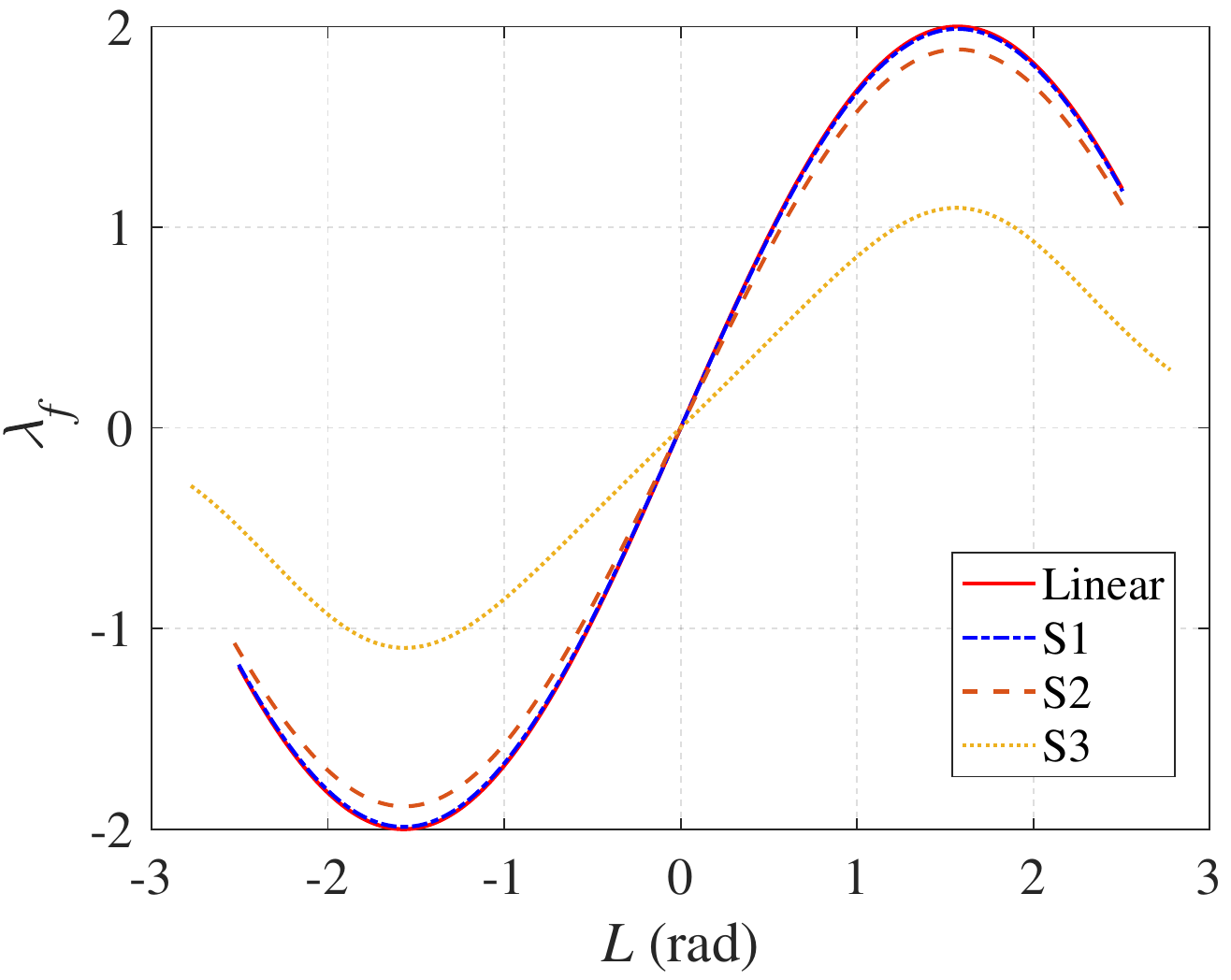} \quad
	\includegraphics[scale = 0.56]{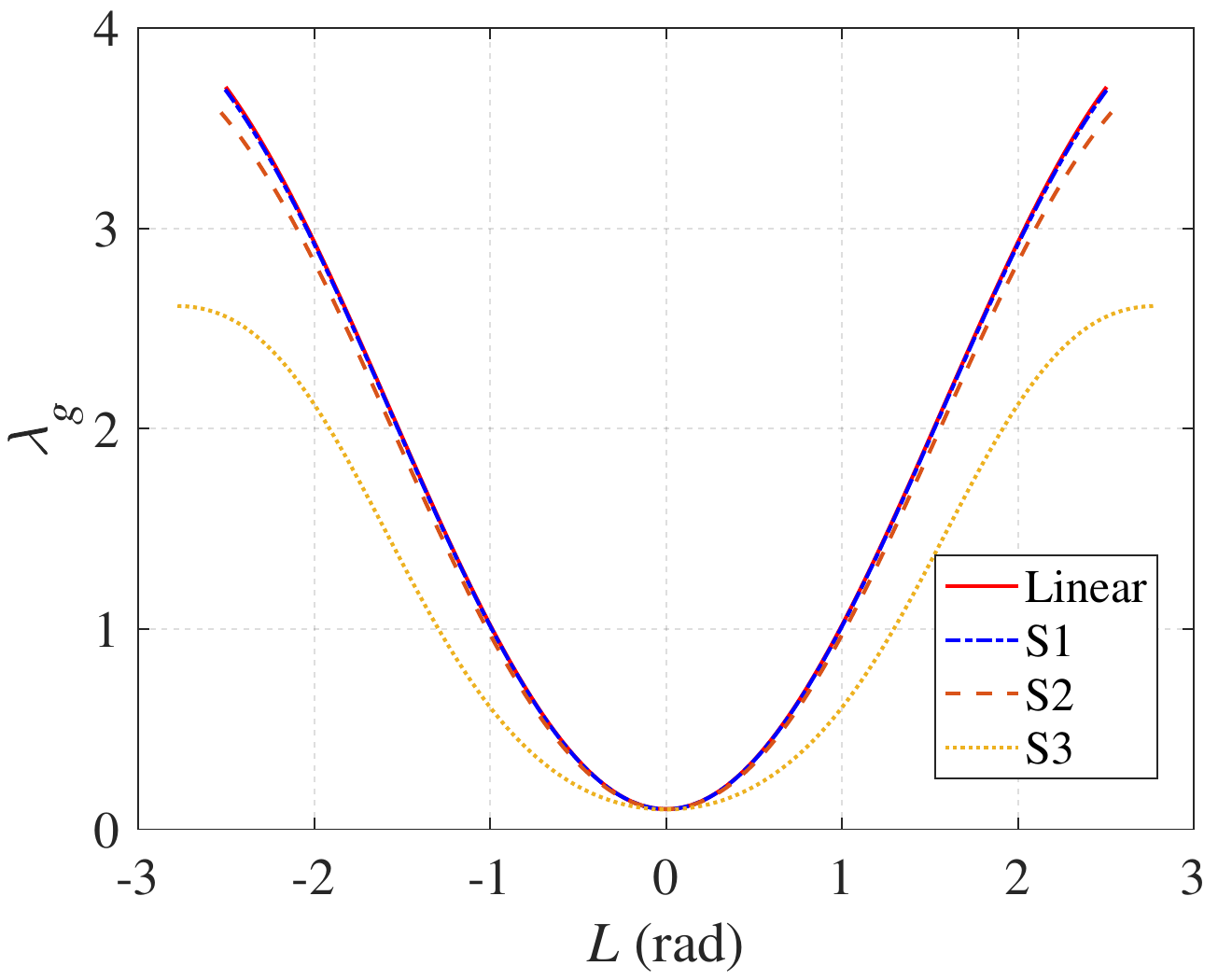}
	\caption{The histories of the control profiles and costates in Case 2.}
	\label{fig9}
\end{figure}
%

%
\begin{figure}[ht!]
	\centering
	\includegraphics[scale = 0.56]{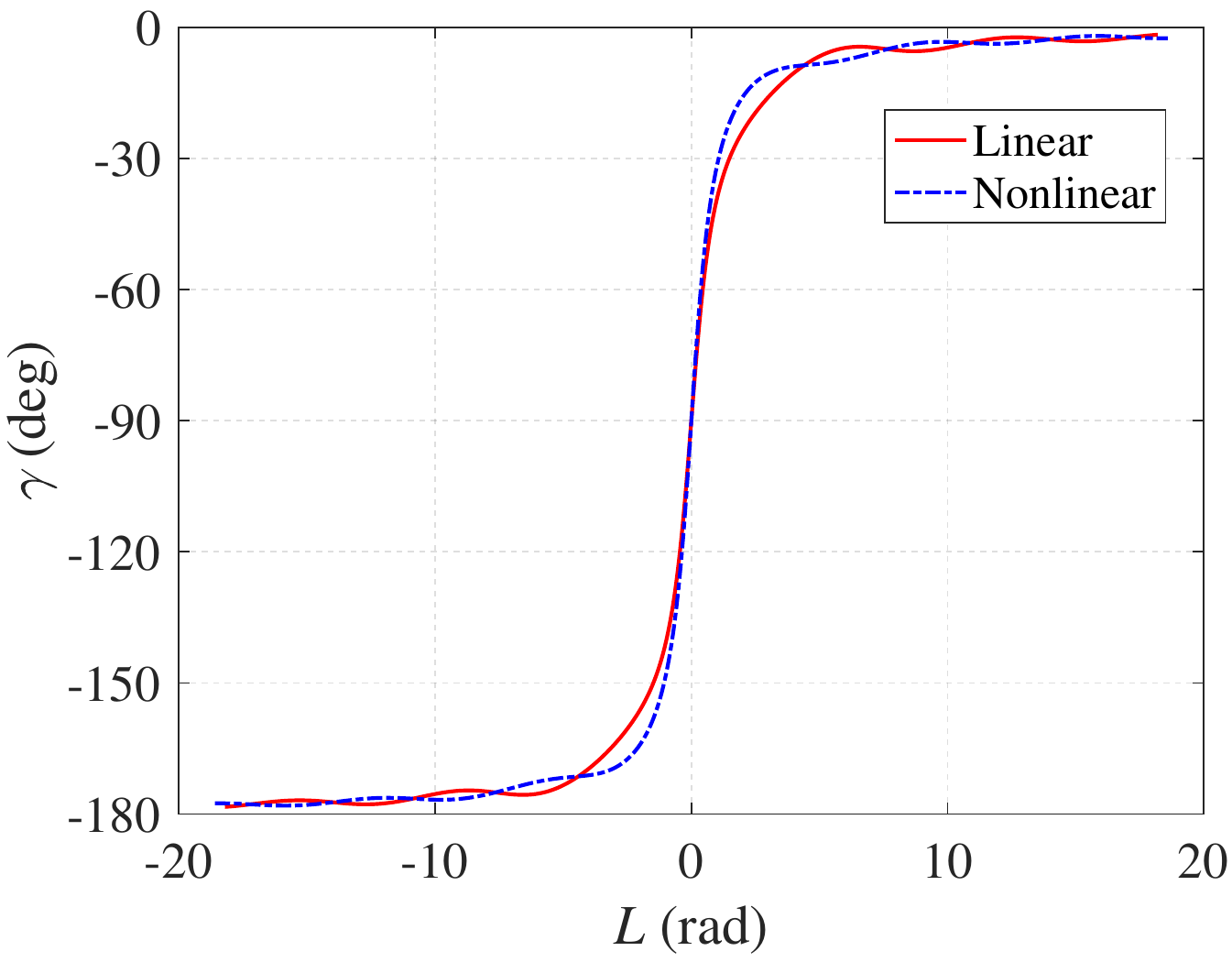} \quad
	\includegraphics[scale = 0.56]{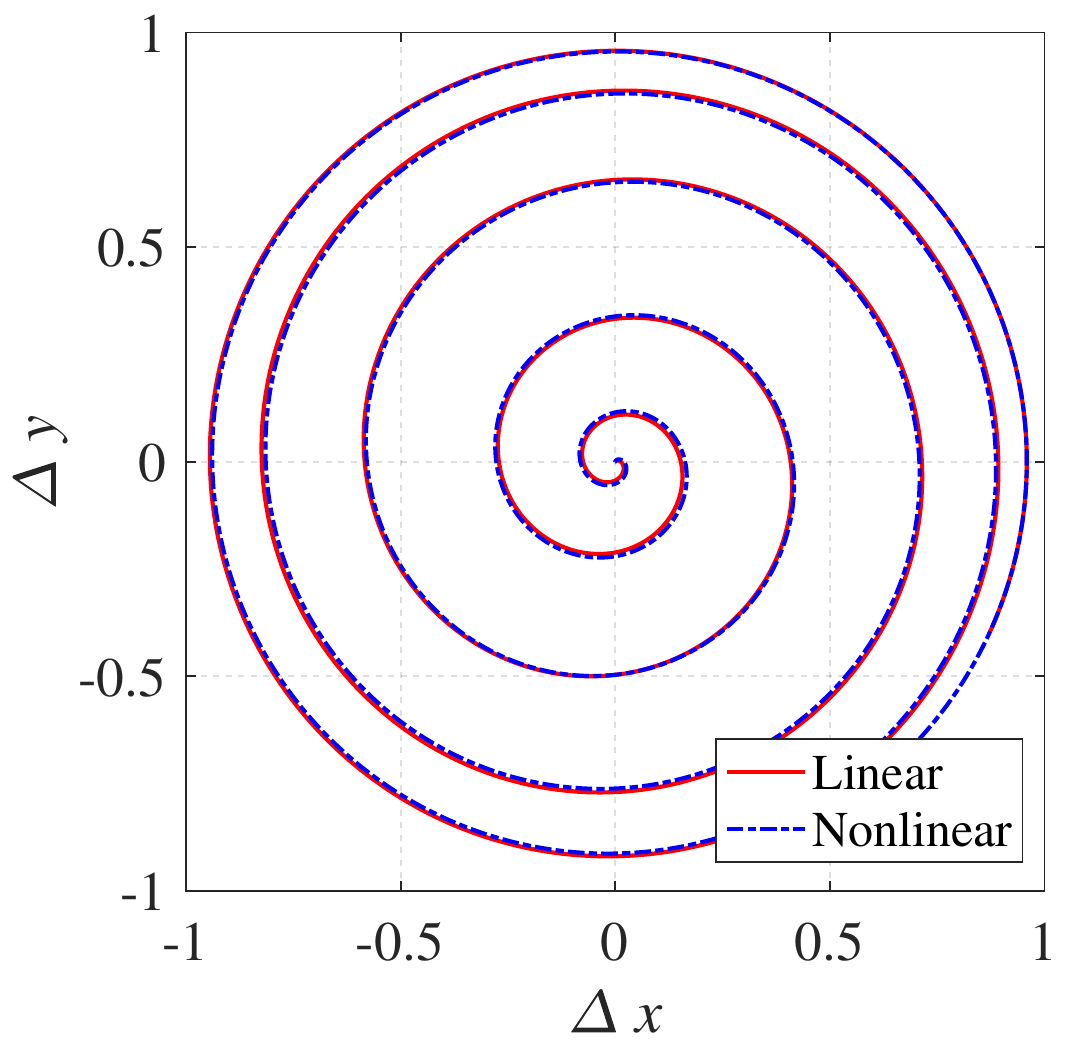} \quad
	\caption{The control profiles and relative trajectories in $xy$ plane in Case 3.}
	\label{fig10}
\end{figure}
%

As for the time-optimal solution, the magnitude of the thrust acceleration takes its maximum value, and the control profile is identified by the orientation angle $\gamma$. The histories of $\gamma$ in each case are shown in Figs.~\ref{fig8}--\ref{fig10}, confirming the bang-bang structure of transversal acceleration in the short- and long-term problems. At the mid-point, the orientation angle is $-90^{\circ}$ in both Case 1 and Case 3, and $90^{\circ}$ in Case 2. All linear and nonlinear solutions seem to satisfy the symmetry property. Figures~\ref{fig8} and~\ref{fig10} present the trajectories of the satellite relative to the target. Note that the relative position at the initial true longitude can be different because of the different $\Delta L$ obtained by the linear and nonlinear solutions. It shows that the overall linear and nonlinear trajectories are quite close in Case 1 and 3. Thus, the proposed linear equations of motion are applicable for the cases with large phase difference. The histories of the costate values are shown in Fig.~\ref{fig9}. The analytical Euler-Lagrange equation~\eqref{Eul1_Solu} accurately approximate the costates in S1 and S2, and it become worse as the thrust acceleration increases.

\subsection{B. Linear and Nonlinear Propellant-Optimal Rephasing Solutions}

Compared with the time-optimal solution, the propellant-optimal solution is more complicated in the control profile, which is determined by the true longitude difference $\Delta L$, the parameter $\eta$, the final time (or phase) difference $\Delta t_f$, the maximum thrust acceleration $a_{\max}$, and the smoothing parameter $\epsilon$. In this subsection, the smoothing parameter $\epsilon$ is set to $0.01$, and a continuation technique introduced in Ref.~\cite{taheri2018generic} is employed to obtain the optimal solution corresponding to smaller smoothing parameter $\epsilon = 1\times 10^{-6}$. The maximum thrust acceleration $a_{\max}$ is set to 0.001, and its effect has been discussed in the time-optimal solutions. Three cases with $\Delta L = 0.5, \, 8, $ and $50$ are tested to compare the linear with nonlinear propellant-optimal solutions, where the parameter $\eta$ takes the values $0.4, \, 0.6$, and $0.8$, respectively. The final time difference $\Delta t_f$ in each case is computed according to Eq.~\eqref{eta} and summarized in Table~\ref{table3}.

%
%

%
{
\begin{table}[!htb]
	\caption{Parameters of the linear and nonlinear propellant-optimal solutions}
	\label{table3}
	\vspace{0.25cm}
	\centering \small
	\begin{tabular}{cccc}
		\hline \hline
		Parameter & Case 1 & Case 2 & Case 3 \\ \hline
		$\Delta L$   & 0.5 &  8 &  50 \\
		$\eta$       & 0.4 &  0.6 & 0.8 \\
		$\Delta t_f$ & $-5.21 \times 10^{-5}$  & $-2.73 \times 10^{-2}$ & $-0.677$\\
		\hline \hline
	\end{tabular}
\end{table}}
%

%
{
\begin{table}[!htb]
	\caption{Comparison of the linear and nonlinear propellant-optimal solutions}
	\label{table4}
	\vspace{0.25cm}
	\centering \small
	\begin{tabular}{ccccc}
		\hline \hline
		Case & Solution & $\V{\lambda}_x$ & $\lambda_t$ & $J \,/\, \left(a_{\max} \, \Delta L\right)$ \\ \hline
		\multirow{3}{*}{1} & Linear & $\left[3.82819,\, -5.05125,\, 0.37921\right]^{\T}$ & 10.20851 &0.61117 \\
		                & Nonlinear & $\left[3.83051,\, -5.05422,\, 0.37959\right]^{\T}$ & 10.21655 & 0.61133 \\
		                & Optimal   & $\left[3.83034, \, -5.05401, \, 0.37956\right]^{\T}$ & 10.21612 & 0.61131 \\
		\multirow{3}{*}{2} & Linear & $\left[0.64131,\, 0.16178,\, -0.03302\right]^{\T}$ & 0.10688 & 0.36119 \\
		                & Nonlinear & $\left[0.64359,\, 0.16434,\, -0.03320\right]^{\T}$ & 0.10776 & 0.36264 \\
		                & Optimal   & $\left[0.64163, \, 0.16384, \, -0.03338\right]^{\T}$ & 0.10743 & 0.36233 \\
		\multirow{3}{*}{3} & Linear & $\left[0.59019,\, 0.00417,\, -0.08094\right]^{\T}$ & 0.01574 & 0.20261 \\
		                & Nonlinear & $\left[0.59382,\, 0.00488,\, -0.08402\right]^{\T}$ & 0.01614 & 0.20485\\
		                & Optimal   & $\left[0.59265, \, 0.00488, \, -0.08439\right]^{\T}$ & 0.01611 & 0.20481 \\
		\hline \hline
	\end{tabular}
\end{table}}
%

The linear solution for each case is first obtained by solving a two-dimensional shooting function~\eqref{shooting2}, in which the initial guess is provided by the linear interpolation. The costate values at the initial true longitude are computed by Eq.~\eqref{Eul1_Solu}. Bearing in mind $\lambda_t = \lambda_0$, the nonlinear solution is then obtained by solving the four-dimensional shooting function~\eqref{nonboundary}. All cases converge when their initial guesses are generated by the linear solution. The costates and performance indexes (or propellant consumptions) of the linear, nonlinear, and optimal solutions are compared in Table~\ref{table4}. The optimal solution corresponding to $\epsilon = 1\times10^{-6}$ is obtained by the continuation technique. In Cases 1 and 3, when the true longitude difference $\Delta L$ is relatively small or large, the propellant consumption is about $J \,/\, \left(a_{\max} \, \Delta L\right) \approx 1 - \eta$, verifying the analytical solutions in Sec.~\uppercase\expandafter{\romannumeral4}.C. In the three cases, the values of costate $\V{\lambda}_x$ and $\lambda_t$ are rather close between the linear and nonlinear solutions. The relative error of the propellant consumption becomes larger as the true longitude $\Delta L$ increases, and its maximum value is about $1.0\%$. Thus, the proposed linear propellant-optimal solutions provide good initial guesses to the nonlinear solutions and accurate approximations to the optimal propellant consumptions. A linear interpolation of contour results presented in Fig.~\ref{fig5}.b is suggested to obtain efficient approximations.

%
\begin{figure}[ht!]
	\centering
	\includegraphics[scale = 0.56]{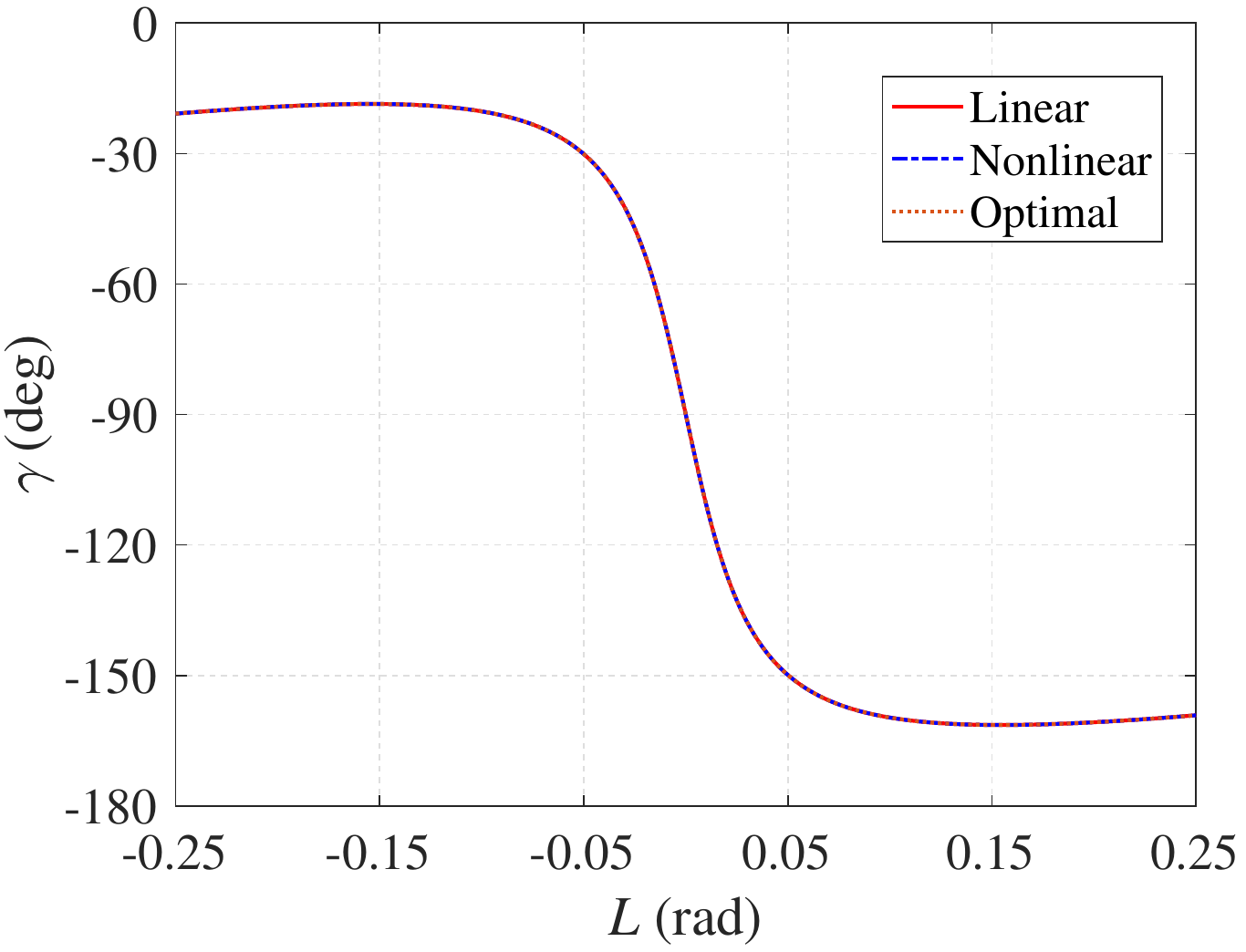} \quad
	\includegraphics[scale = 0.56]{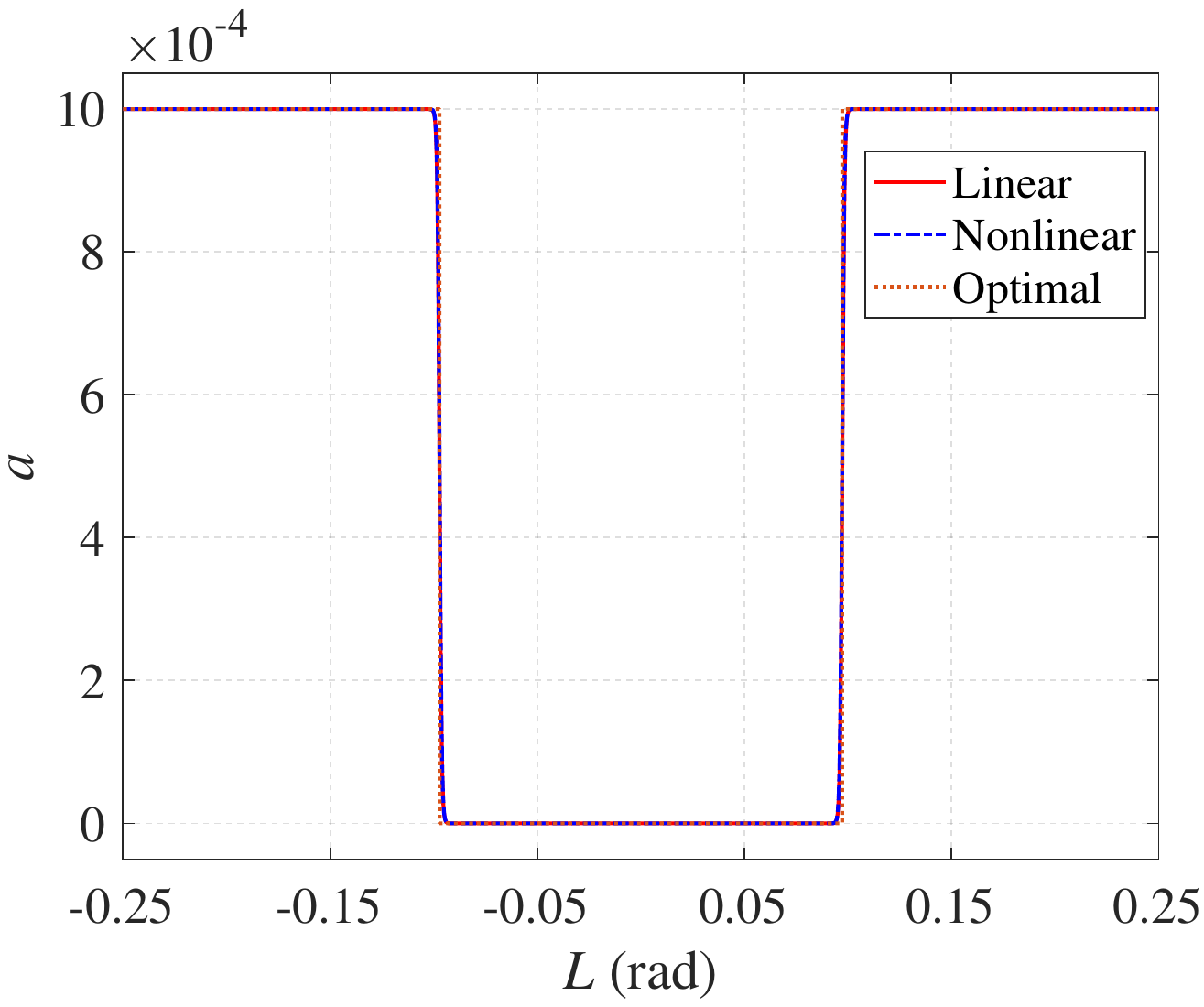} \\
	\;\includegraphics[scale = 0.56]{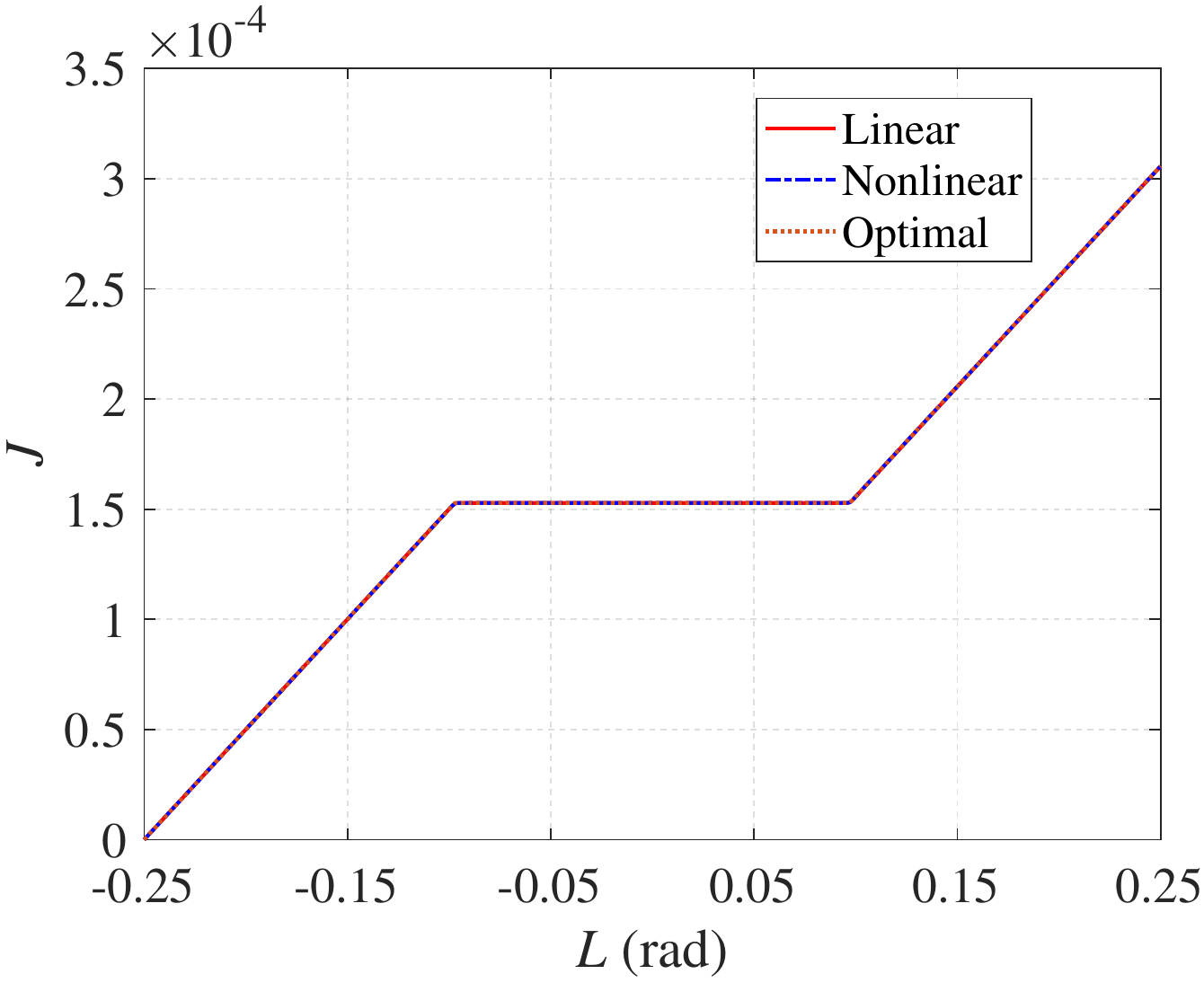} \quad\;\qquad
	\includegraphics[scale = 0.56]{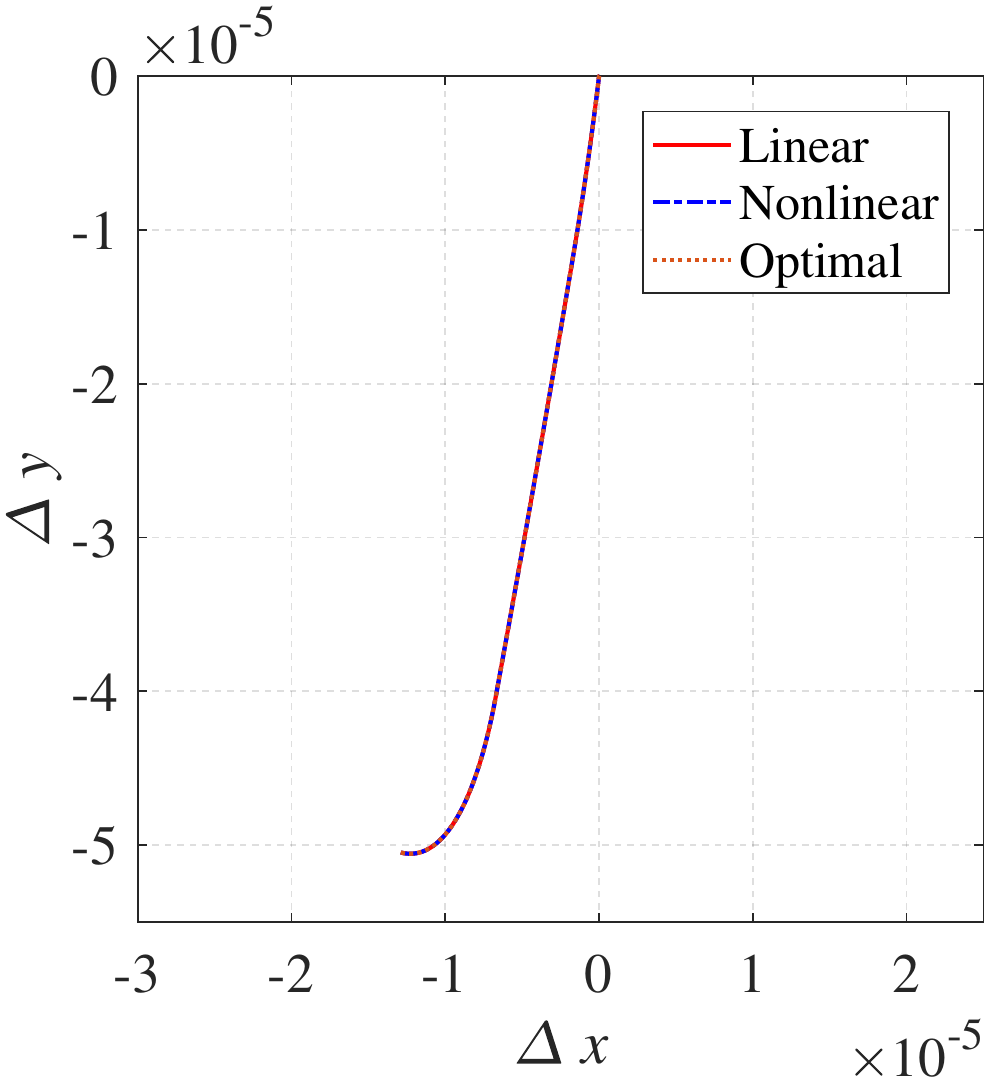} \quad\,
	\caption{The control profiles, propellant consumptions, and relative trajectories in Case 1.}
	\label{fig11}
\end{figure}
%

%
\begin{figure}[ht!]
	\centering
	\includegraphics[scale = 0.56]{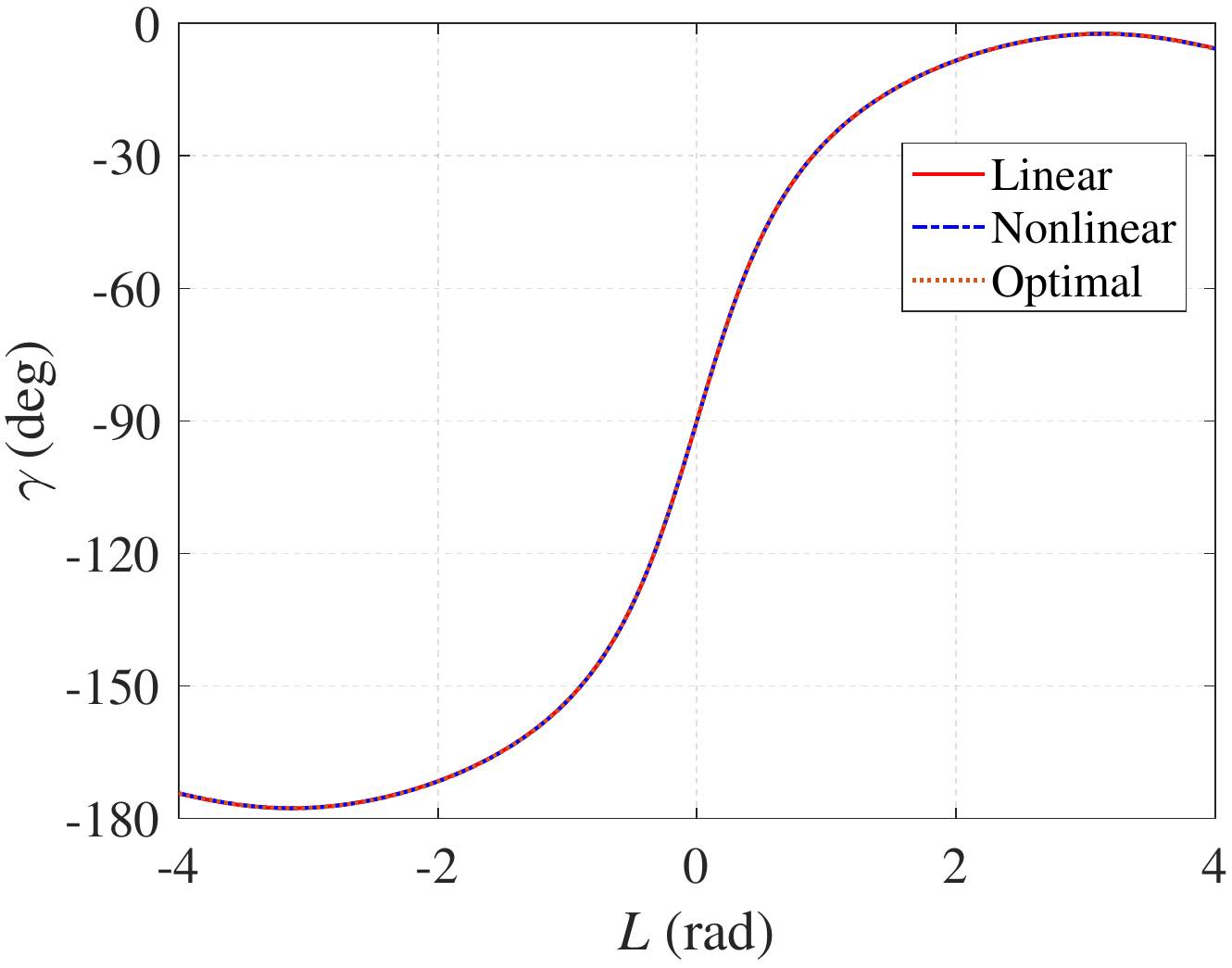} \quad\;
	\includegraphics[scale = 0.56]{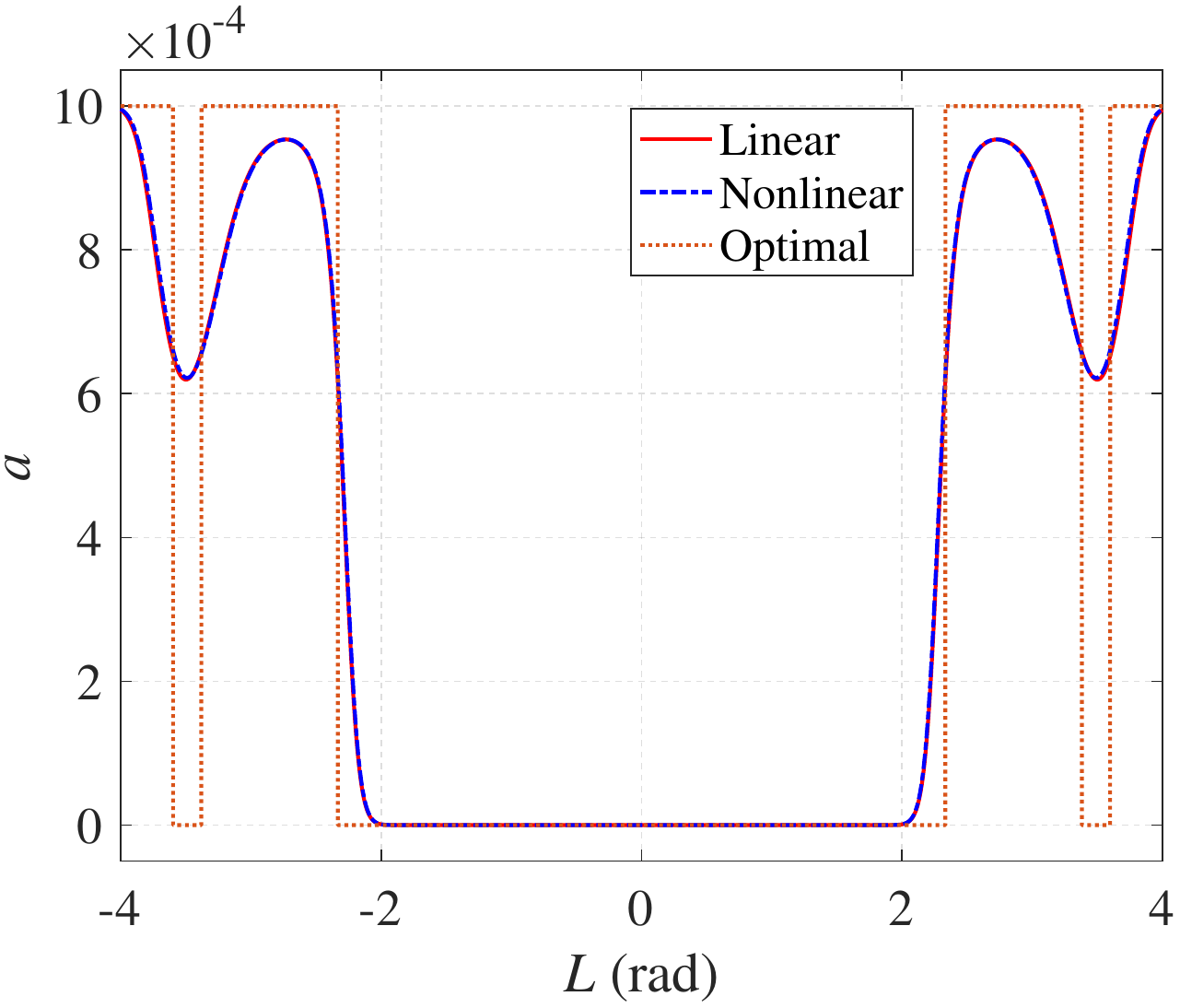} \\
	\;\includegraphics[scale = 0.56]{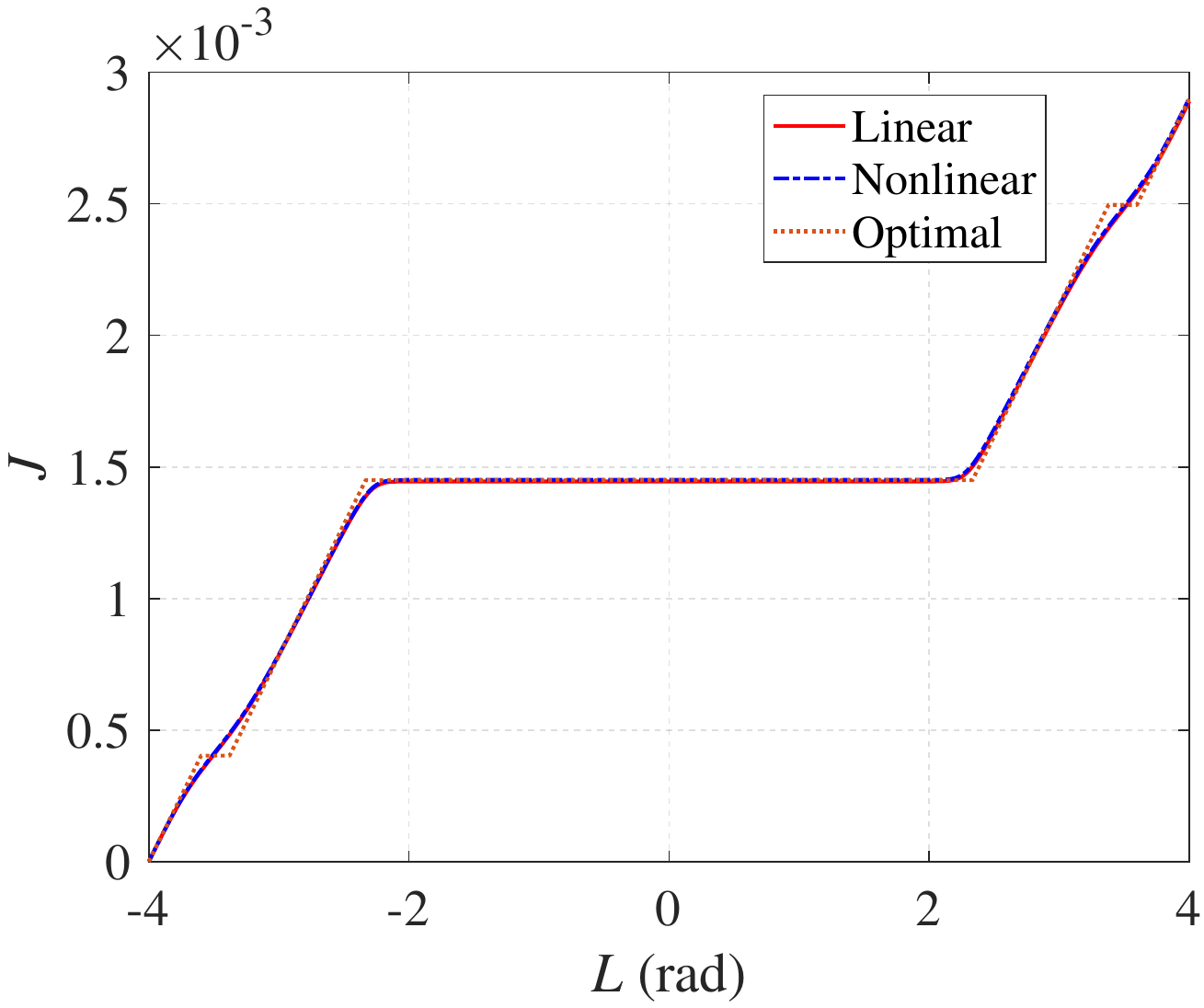} \quad
	\includegraphics[scale = 0.56]{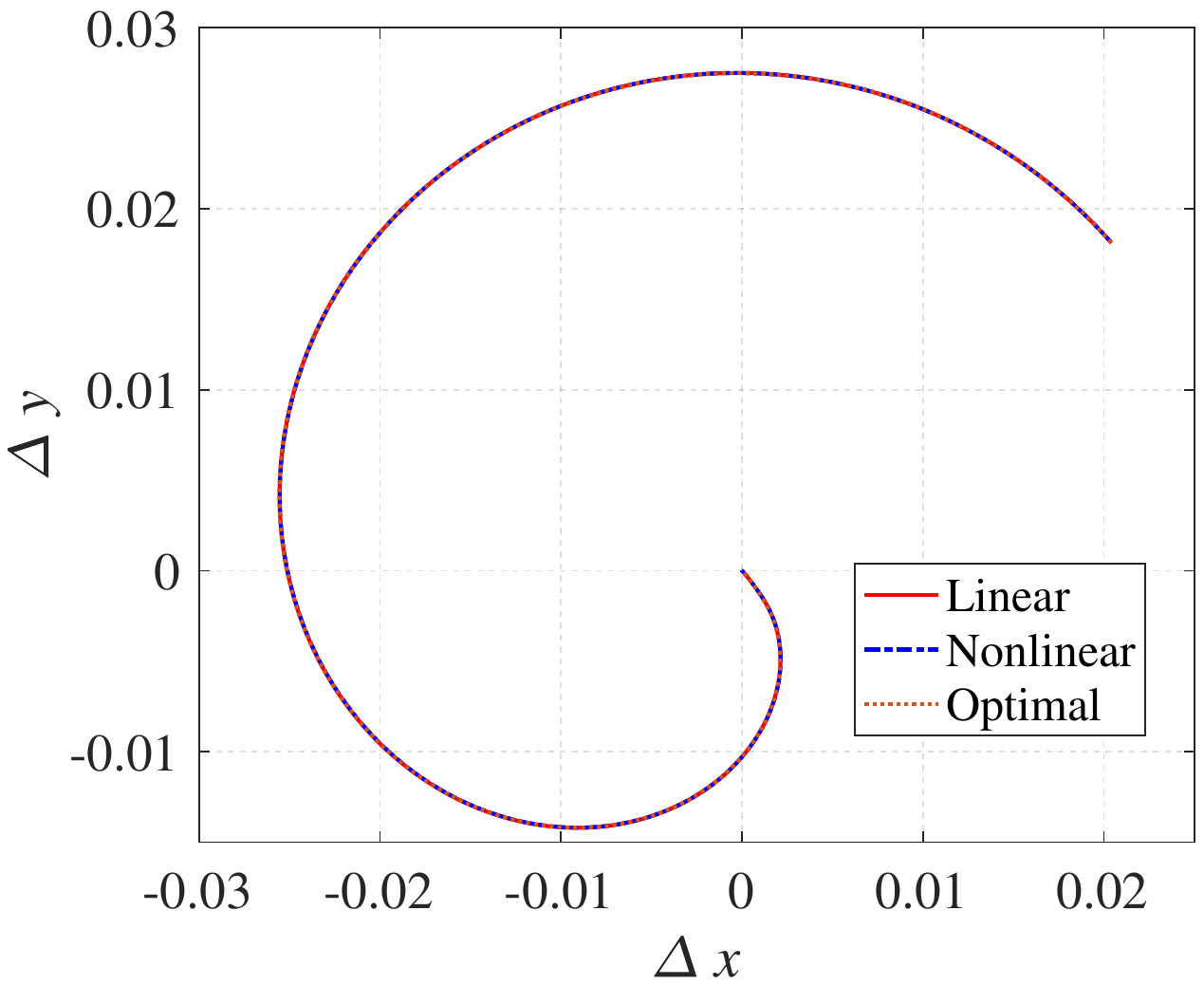}
	\caption{The control profiles, propellant consumptions, and relative trajectories in Case 2.}
	\label{fig12}
\end{figure}
%

%
\begin{figure}[ht!]
	\centering
	\includegraphics[scale = 0.56]{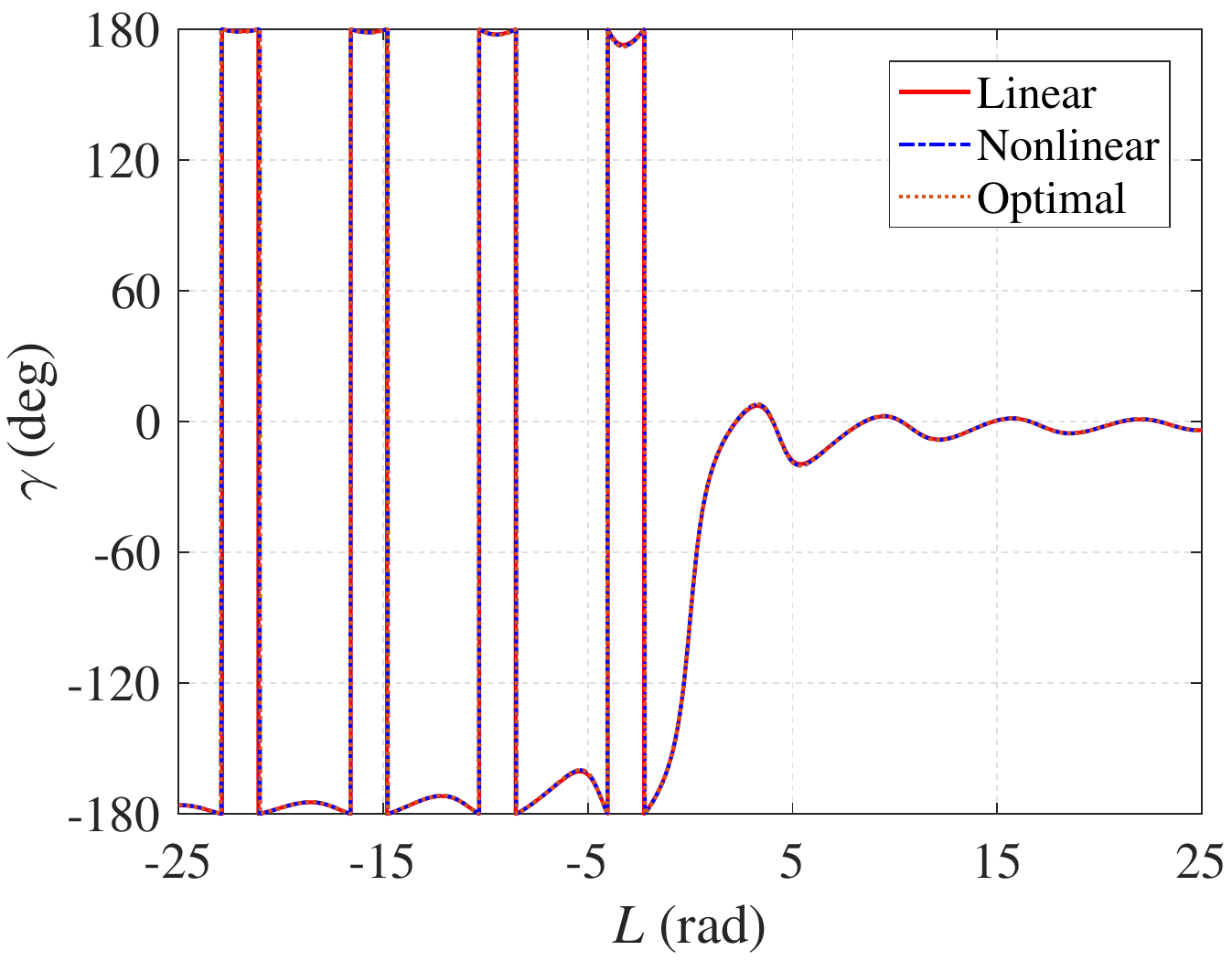} \quad
	\includegraphics[scale = 0.56]{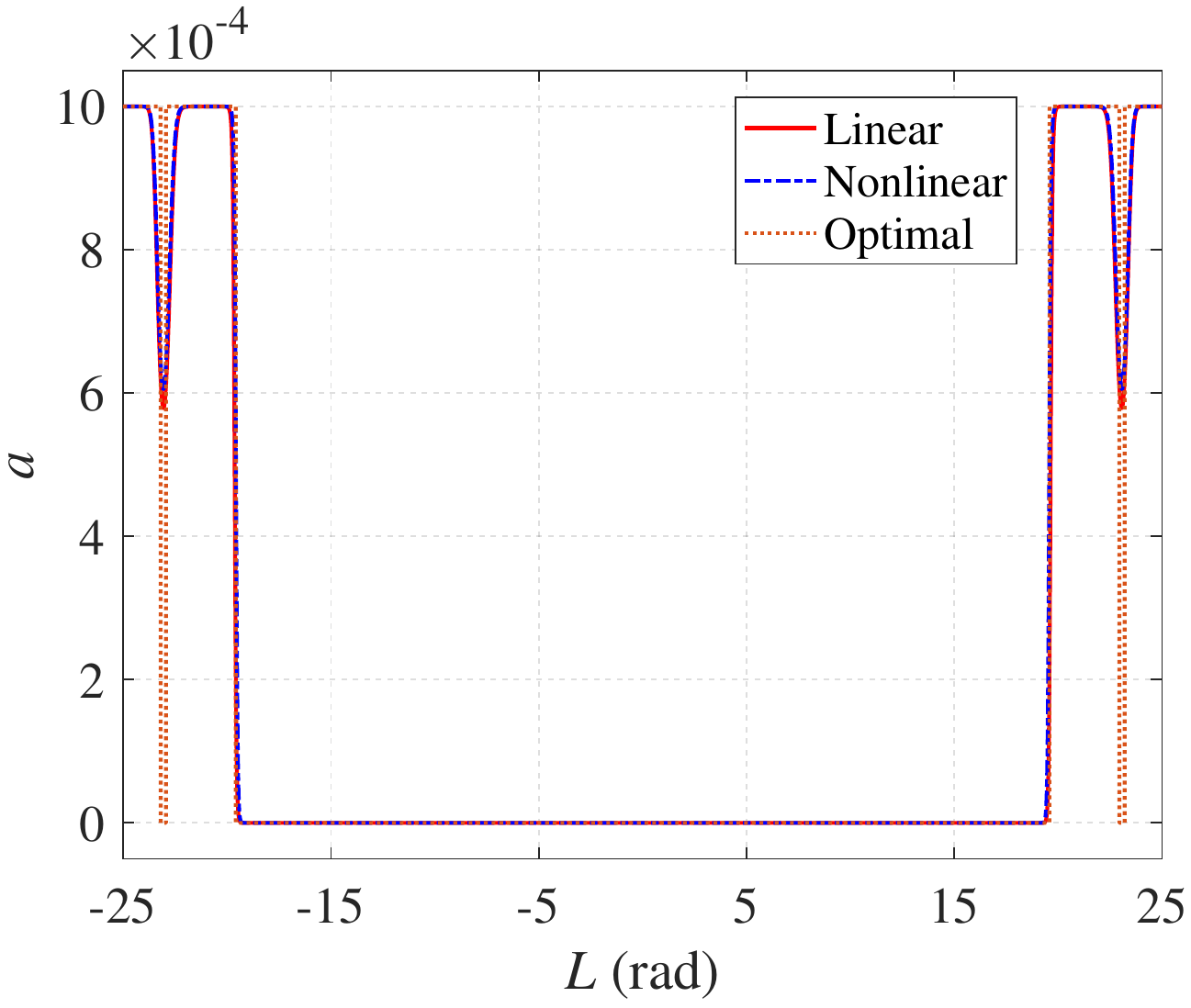} \\
	\;\;\includegraphics[scale = 0.57]{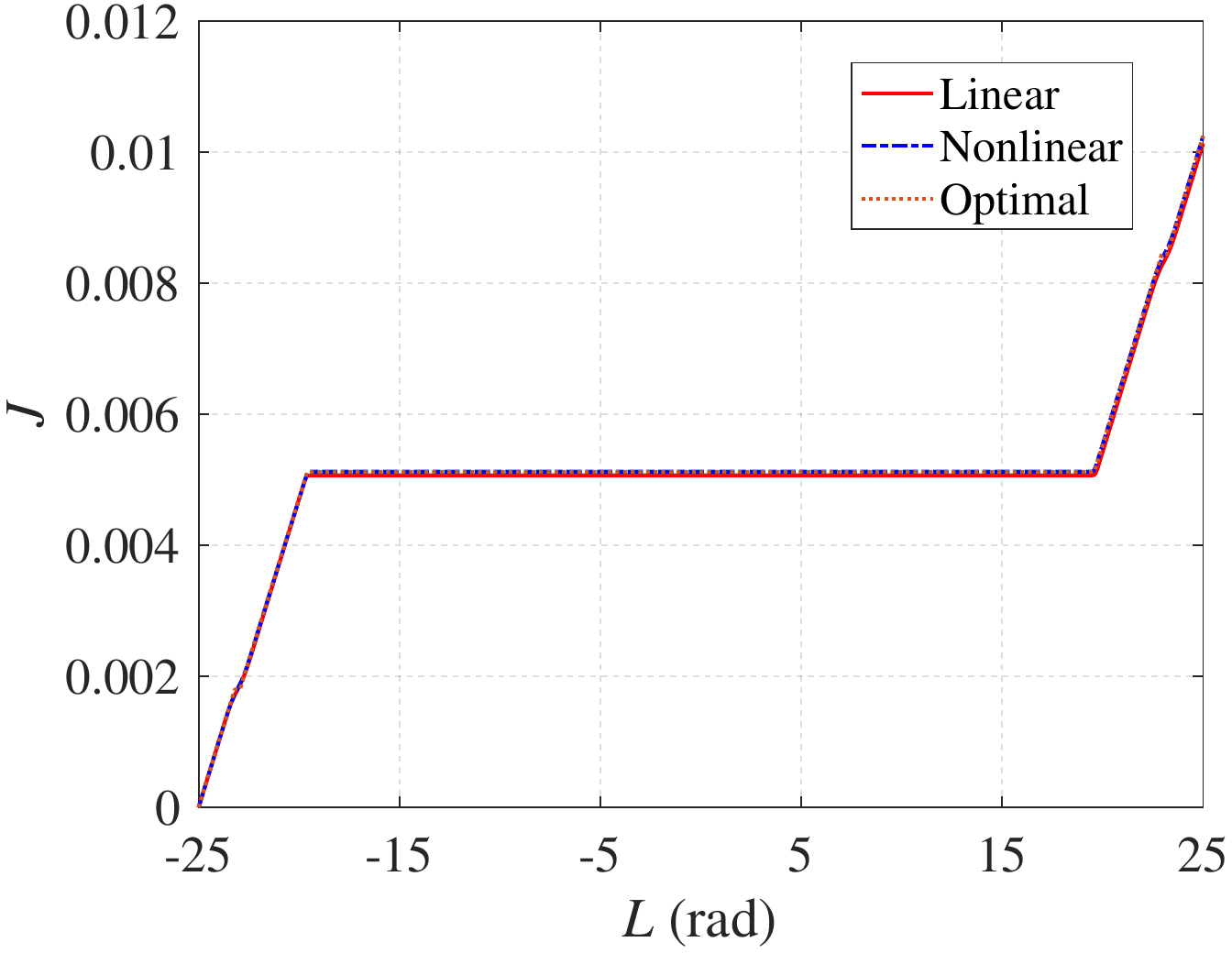} \qquad\quad
	\includegraphics[scale = 0.57]{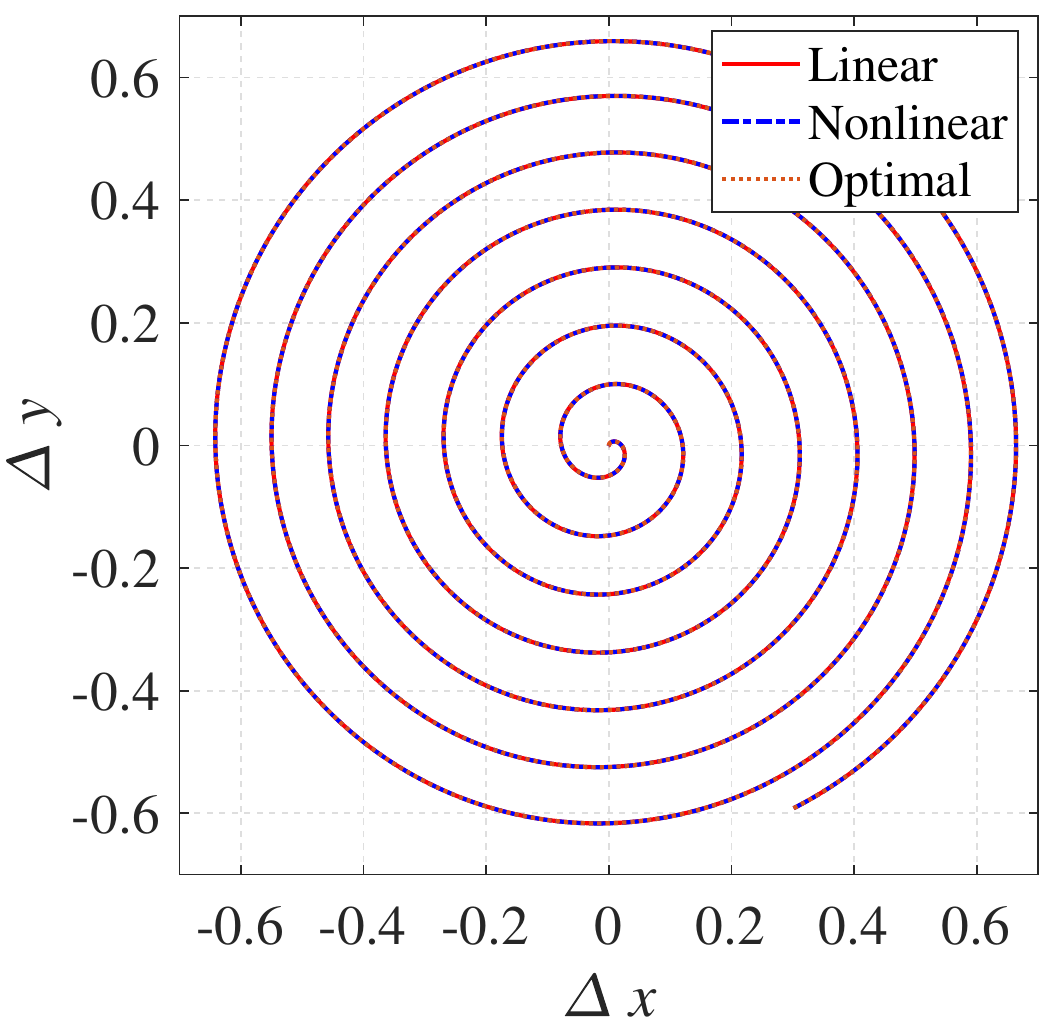} \quad\;\quad
	\caption{The control profiles, propellant consumptions, and relative trajectories in Case 3.}
	\label{fig13}
\end{figure}
%

The results of the optimal control profiles, propellant consumptions, and trajectories of the three cases are depicted in Figs.~\ref{fig11}--\ref{fig13}. The control profiles are characterized by the thrust orientation angle $\gamma$ and the smoothed magnitude $a^{\star}$. In the simulation, the orientation angle is computed independently of the magnitude, and it is meaningless at the coast arcs where $a^{\star} = 0$. The optimal magnitudes of the linear and nonlinear solutions with $\epsilon = 0.01$ are smoothed, and the bang-bang control is obtained by solving the problem with $\epsilon = 1\times10^{-6}$ by the continuation technique~\cite{taheri2018generic}. As shown in Fig.~\ref{fig11}, the linear and nonlinear results are almost the same. Although the phase difference is quite small in this case, these results demonstrate the closeness between the uses of linear and nonlinear equations to solve the low-thrust rephasing problem. In Figs.~\ref{fig12} and~\ref{fig13}, the linear and nonlinear solutions look close to each other. The optimal linear and nonlinear solutions have similar thrust orientations $\gamma$ and trajectories, and the optimal thrust magnitudes are of bang-bang control structure. The number of burning arcs is two in Case 1, in good agreement with the analytical solution to the short-term problem. The number of burning arcs is four in Case 3, where the two short coast arcs are neglected by the analytical solution to the long-term problem. In general, the proposed analytical solutions for the short- and long-term rephasing are near-optimal.

\section{VI. Conclusions}

This work presented an atlas of the time- and propellant-optimal low-thrust rephasing solutions in the circular orbit, which provide good initial guesses to solve the nonlinear problem and approximated performance indexes for the preliminary mission design. These solutions are obtained by numerically solving two reduced two-dimensional shooting functions based on a set of linearized and scaled equations of motion. The time-optimal rephasing solutions are dependent on only one key parameter and can be estimated by two piecewise functions. The number of the key parameters for propellant-optimal solutions is two, and the solution space is described by some contour maps in combination with the linear interpolation technique. The switching functions are investigated to explain why the costate values mutate with the boundary constraints. The number of the burning arcs is found to be two or four for all propellant-optimal solutions. In addition, the analytical solutions to the relatively short- and long-term problems are derived and in good agreement with the numerical results.

Numerical test shows that the proposed time- and propellant-optimal solutions can be efficiently estimated and numerically obtained by the shooting algorithm in several iterations. The costates values, optimal control profiles, and trajectories are close to those obtained with the nonlinear dynamics, respectively, assuming that the low-thrust acceleration is high-order smaller than the gravitational acceleration. The solutions with the smoothing parameter $\epsilon = 0.01$ are near-optimal compared with the propellant-optimal solutions, and the bang-bang control can be achieved by the traditional continuation technique.

\section{Acknowledgment}
This work was supported by the National Natural Science Foundation of China (Grants. U21B2050 and 12022214).


\bibliographystyle{aiaa_doi}               
\bibliography{biblio_di}

\end{spacing}
\end{document}